\documentclass[a4paper,10pt,twocolumn]{article}
\usepackage[utf8]{inputenc}

\usepackage[switch]{lineno}

\usepackage{hyperref}
\usepackage{bm}
\usepackage{authblk}
\usepackage{siunitx}
\usepackage{graphicx}

\usepackage[a4paper,left=1.5cm,right=1.5cm,top=2cm,bottom=2cm]{geometry}

\usepackage{nomencl}
\usepackage{ifthen}
\renewcommand{\nomgroup}[1]{%
\ifthenelse{\equal{#1}{C}}{\item[\textbf{Constants}]}{%
\ifthenelse{\equal{#1}{G}}{\item[\textbf{Greek Letters}]}{%
\ifthenelse{\equal{#1}{N}}{\item[\textbf{Dimensionless Numbers}]}{%
\ifthenelse{\equal{#1}{S}}{\item[\textbf{Subscripts}]}{}}}}
}
\makenomenclature

\nomenclature[G]{${\rho}$}{ Density [$\si{kg/m^3}$]}
\nomenclature[G]{${\beta}$}{ Fabric Properties [$\si{-}$]}
\nomenclature[G]{${\chi}$}{ Water Content [$\si{m^3}$]}
\nomenclature[G]{${\epsilon}$}{ Effectiveness [$\si{-}$]}
\nomenclature[G]{${\nu}$}{ Momentum diffusivity [$\si{m^2/s}$]}
\nomenclature[G]{${\alpha}$}{ Thermal diffusivity [$\si{m^2/s}$]}
\nomenclature[G]{${\lambda}$}{ Thermal Conductivity [$\si{W/mK}$]}
\nomenclature[G]{${\tau}$}{ Non-dimensional Residence Time [$\si{-}$]}

\nomenclature[N]{Gu}{ Gukhman number }
\nomenclature[N]{Le}{ Lewis number }
\nomenclature[N]{Nu}{ Nusselt number }
\nomenclature[N]{Pr}{ Prandtl number }
\nomenclature[N]{Sc}{ Schmidt number }
\nomenclature[N]{Sh}{ Sherwood number}
\nomenclature[N]{Sh\textsubscript{a}}{ Area-Sherwood number }
\nomenclature[N]{Re}{ Reynolds number }
\nomenclature[N]{Le'}{ Modified Lewis number }

\nomenclature{$A$}{ Area [$\si{m^2}$]}
\nomenclature{$a$}{ Activity Coefficient [$\si{-}$]}
\nomenclature{$C_p$}{Specific heat [$\si{J/kg K}$]}
\nomenclature{$C(t)$}{ Instantaneous moisture content [\%]}
\nomenclature{$d$}{ Diameter [$\si{m}$]}
\nomenclature{$D$}{ Diffusion coefficient [$\si{m^2/s}$]}
\nomenclature{$g$}{ Gravity [$\si{m/s^2}$]}
\nomenclature{$h$}{ Heat transfer coefficient [$\si{W/m^2K}$]}
\nomenclature{$h_{fg}$}{ Latent heat of evaporation  [$\si{kJ/kg}$]}
\nomenclature{$k$}{ Mass transfer coefficient [$\si{m/s}$]}
\nomenclature{$kA$}{ Area–Mass Transfer Coefficient [$\si{m^3/s}$]}
\nomenclature{$k_B$}{Boltzmann Constant [$\si{J/K}$]}
\nomenclature{$L$}{ Length [$\si{m}$]}
\nomenclature{$m$}{ Weight [$\si{kg}$]}
\nomenclature{$\dot{m}$}{ Mass Flow rate [$\si{kg/s}$]}
\nomenclature{$N$}{ Drum Angular Speed [$\si{rpm}$]}
\nomenclature{$NTU$}{ Number of transfer units [$\si{-}$]}
\nomenclature{$P$}{ Pressure [$\si{Pa}$]}
\nomenclature{$PO$}{ Power [$\si{W}$]}
\nomenclature{$\dot{Q}$}{ Heat transfer rate [$\si{W}$]}
\nomenclature{$RH$}{ Relative humidity [{\%}]}
\nomenclature{$T$}{ Temperature [$\si{K}$]}
\nomenclature{$t$}{ Time [$\si{sec}$]}
\nomenclature{$V$}{ Volume [$\si{m^3}$]}
\nomenclature{$\dot{V}$}{ Flow rate [$\si{m^3/s}$]}
\nomenclature{$w$}{ Humidity ratio [$\si{kg/kg}$]}
	
\nomenclature[S]{$a$ }{ Air}
\nomenclature[S]{$amb$ }{ Ambient}
\nomenclature[S]{$d$ }{ Drum}
\nomenclature[S]{$clo$ }{ Clothes}
\nomenclature[S]{$evap$ }{ Evaporated}
\nomenclature[S]{$ext$ }{ External}
\nomenclature[S]{$f$ }{ Flow}
\nomenclature[S]{$g$ }{ Dry bulb}
\nomenclature[S]{$gw$ }{ Wet bulb}
\nomenclature[S]{$h$ }{ Heater}
\nomenclature[S]{$LH$ }{ Latent Heat}
\nomenclature[S]{$in$ }{ Inlet}
\nomenclature[S]{$mois$}{ Moisture}
\nomenclature[S]{$out$ }{ Outlet}
\nomenclature[S]{$sat$ }{ Saturation}
\nomenclature[S]{$surf$ }{ Surface}
\nomenclature[S]{$tot$ }{ Total}
\nomenclature[S]{$t$ }{ Textile}
\nomenclature[S]{$th$ }{ Thickness}
\nomenclature[S]{$v$ }{ Vapor}
\nomenclature[S]{$w$ }{ Water}

\usepackage[authoryear]{natbib}

\begin{document}
\let\WriteBookmarks\relax
\def\floatpagepagefraction{1}
\def\textpagefraction{.001}

\title{Applicable Methodologies for the Mass Transfer Phenomenon in Tumble Dryers: A Review}
\author[1]{Sajad~Salavatidezfouli\footnote{ssalavat@sissa.it}}
\author[1]{Arash~Hajisharifi\footnote{ahajisha@sissa.it}}
\author[1]{Michele~Girfoglio\footnote{mgirfogl@sissa.it}}
\author[1]{Giovanni~Stabile\footnote{gstabile@sissa.it}}
\author[1]{Gianluigi~Rozza\footnote{grozza@sissa.it}}

\affil[1]{Mathematics Area, mathLab, SISSA, via Bonomea 265, I-34136 Trieste, Italy}

\date{} 

\twocolumn[
  \begin{@twocolumnfalse}
    \maketitle
	   \begin{abstract}
        Tumble dryers offer a fast and convenient way of drying textiles independent of weather conditions and therefore are frequently used in ordinary households. However, artificial drying of textiles consumes considerable amounts of energy, approximately  8.2\% of the residential electricity consumption is for drying of textiles in northern European countries \citep{cranston2019efficient}. Several authors have investigated the aspects of the clothes drying cycle with experimental and numerical methods to understand and improve the process. The first turning point study on understanding the physics of evaporation for tumble dryers was presented by \cite{lambert1991modeling} in the early 90s. With the aid of Chilton–Colburn analogy, they introduced the concept of area-mass transfer coefficient to address \textit{evaporation rate}. The evaporation rate is considered to be the main system parameter for dryers with which other performance parameters including drying time, effectiveness, moisture content and efficiency can be estimated. Afterwards, several experimental or numerical studies were published based on this concept, and furthermore, the model was then developed into 0-dimensional \citep{deans2001modelling} and 1-dimensional \citep{wei2017mathematical} to gain more accuracy. 

        More recent literature focused on utilizing dimensional analysis or image processing techniques to correlate drying indices with system parameters. However, the validity  of these regressed models is machine-specific, and hence, cannot be generalized yet. All the previous models for estimating the evaporation rate in tumble dryers are discussed. The review of the related literature showed that all of the previous models for the prediction of the evaporation rate in the clothes dryers have some limitations in terms of accuracy and applicability.
	
    	\vspace{0.5cm}
		\textbf{Highlights:}
		\begin{itemize}
		  \item of washer dryers as the second main domestic consumer of electrical energy,
		
		  \item The experimental results are machine-specific and are hard to be generalized,
		
		  \item Area-mass transfer as the key parameter in evaluating the evaporation rate,
		  \item More complex models have recently been presented by several authors.
    
		\end{itemize}

		\end{abstract}
  \end{@twocolumnfalse}

]
\maketitle
\printnomenclature

\section{Introduction}
The perpetual campaign over the optimal utilization of limited energy resources and the achievement of climate change goals has prompted concentrations on technologies with high energy consumption such as drying. Drying, also known as dehydration or dewatering, is the process of removing moisture from material down to a specific value while maintaining the quality and structure of the product. In traditional methods, the drying took place outdoors under windy and/or sunny conditions. However, in modern society, due to the lack of space and demand for a faster drying time, various artificial methods have been introduced \citep{stawreberg2012energy}. Accordingly, drying systems are broadly used in domestic and industrial applications such as food processing, textile, paper, wood, ceramics, minerals, pharmaceutical, HVAC and ventilation and biotechnology \citep{zahed1995modelling, khaldi2017improving, alishah2018solar, oueslati2018thermal, lowrey2020rotary, salhi2021numerical, wang2021experimental}. 

Drying is roughly responsible for up to 25\% of energy consumption in process industries and is counted as one of the most energy-intensive operations \citep{kudra2004energy,kemp2011fundamentals}. The first utilization of industrial drying systems goes back to chemical engineering a few decades ago \citep{kerkhof2002drying}. This short lifetime along with the high application of the drying process in the industry has not been able to meet today's requirements in terms of the worldwide energy crisis, minimum lifecycle costs and environmentally-friendly technologies \citep{defraeye2014advanced}. 

\begin{figure}
\centering  \vspace*{0.25cm} \hspace*{0cm}
\includegraphics[scale=0.4]{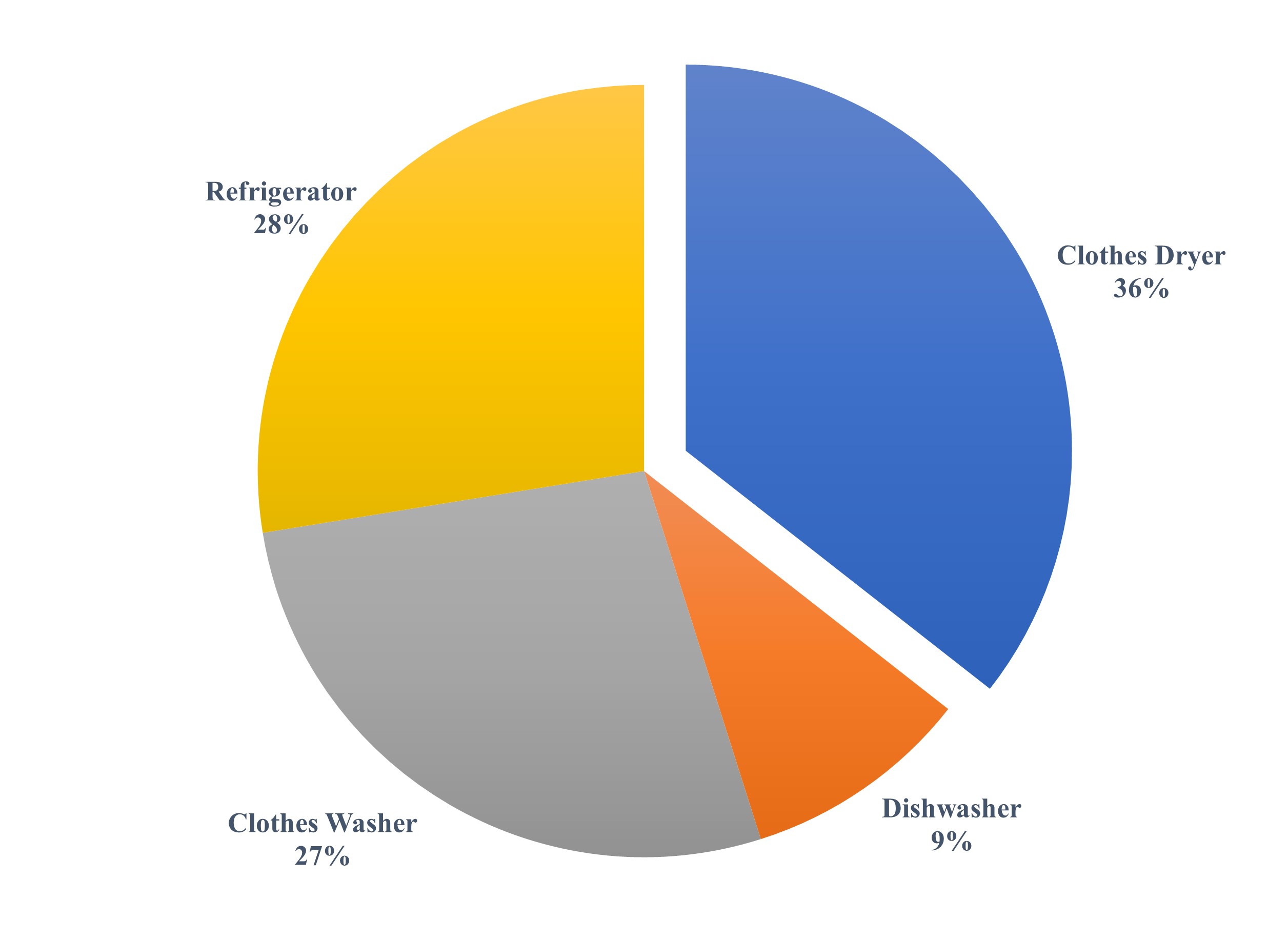}
\caption{Energy consumption of standard home appliances}
\label{EneConsmp}
\end{figure}

One application of drying technology is the drying of porous fabrics such as clothes. The use of domestic tumble clothes dryers has been growing continuously for decades \citep{ng2008new}. Whilst, drying is one of the most energy-consuming processes \citep{minea2015overview}. According to a Department of Energy chart (Fig. \ref{EneConsmp}), the clothes dryer consumes the second largest amount of electricity among household appliances, after the water heater \citep{EnergConsump}. This is quite a bit of electricity use for an appliance that is only run once or twice a week, as compared to an appliance that is running constantly. On this basis, clothes dryers are now considered to be inseparable appliances in ordinary households as they are reliable and fast drying systems \citep{gataric2021influence}. More than 40\% and 80\% of households in Germany and the USA, respectively, are equipped with clothes dryers \citep{EnergyRECS} and this market is steadily increasing \citep{ReviewStudyTumbleDriers}. But, the only disadvantage is their large electrical energy consumption \citep{stawreberg2010modelling}. Estimations show that more than 3\% of residential energy usage in the USA stems from clothes dryers \citep{bassily2003performance} while this number is about 8\% in the world \citep{bansal2010novel}. Tumble dryers are used by people and are gradually growing at an annual rate of 10.9\% due to the changes in consumption structure and living environment (living space and environmental degradation), the increasing health awareness (sun-drying kills only 15\% of the bacteria) and the development of drying-care machinery \citep{MollyPhDThes}. The drying process of fabric is energy-intensive and slow, which takes a long time to complete (up to 3 h, depending on the drying conditions or the dryer model). In addition, energy crisis and environment-friendly have become the focus of citizens’ attention. Therefore, dryer manufacturers have a significant interest in improving the energy efficiency of dehydration to maintain their sale share, competitiveness, and sustainability \citep{wei2017mathematical}.

\begin{figure}
\centering  \vspace*{0.5cm} \hspace*{0cm}
\includegraphics[scale=0.3]{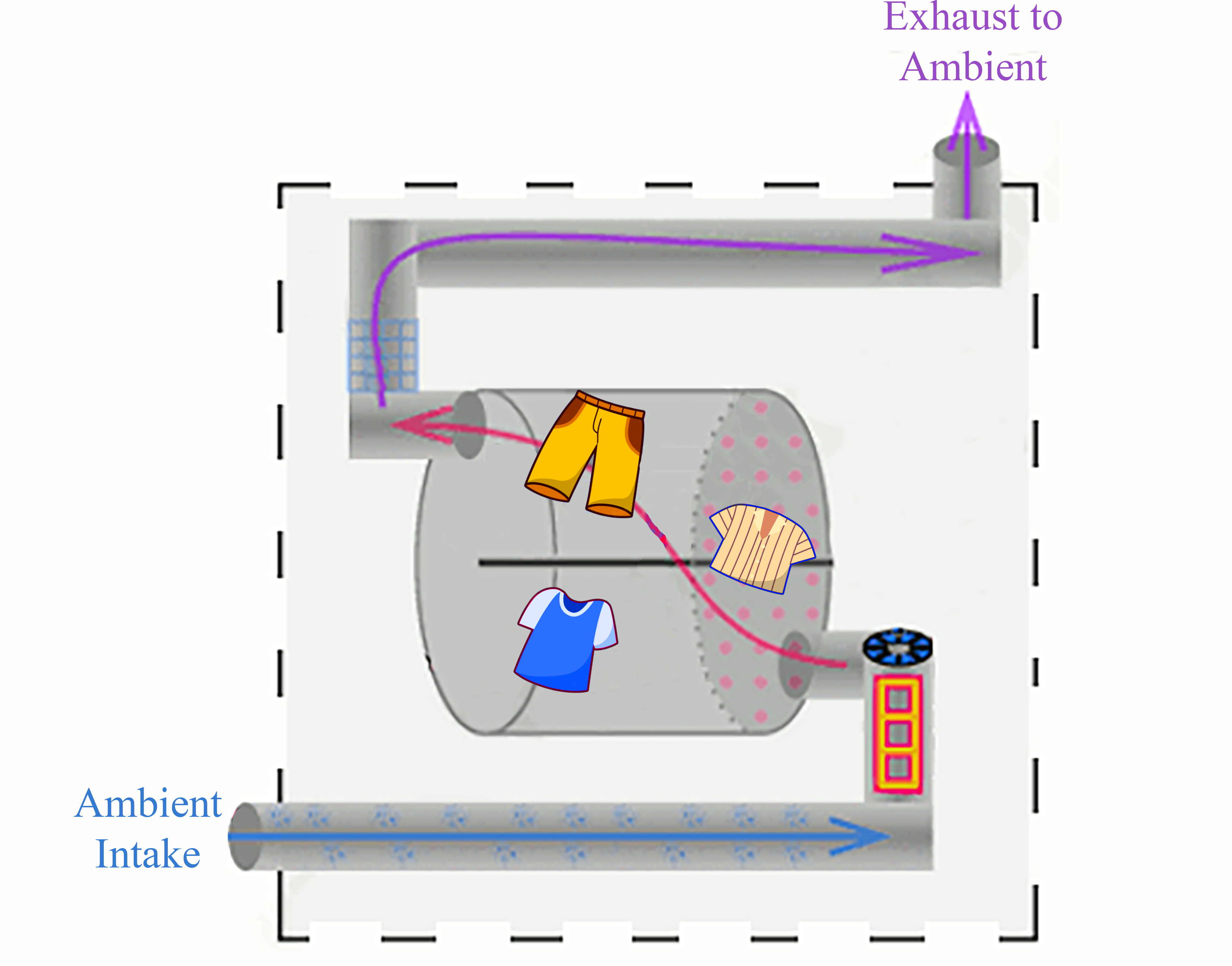}
\caption{Schematic of drying cycle of the vented dryers.}
\label{vented1}
\end{figure}

\begin{figure}
\centering  \vspace*{0.5cm} \hspace*{0cm}
\includegraphics[scale=0.6]{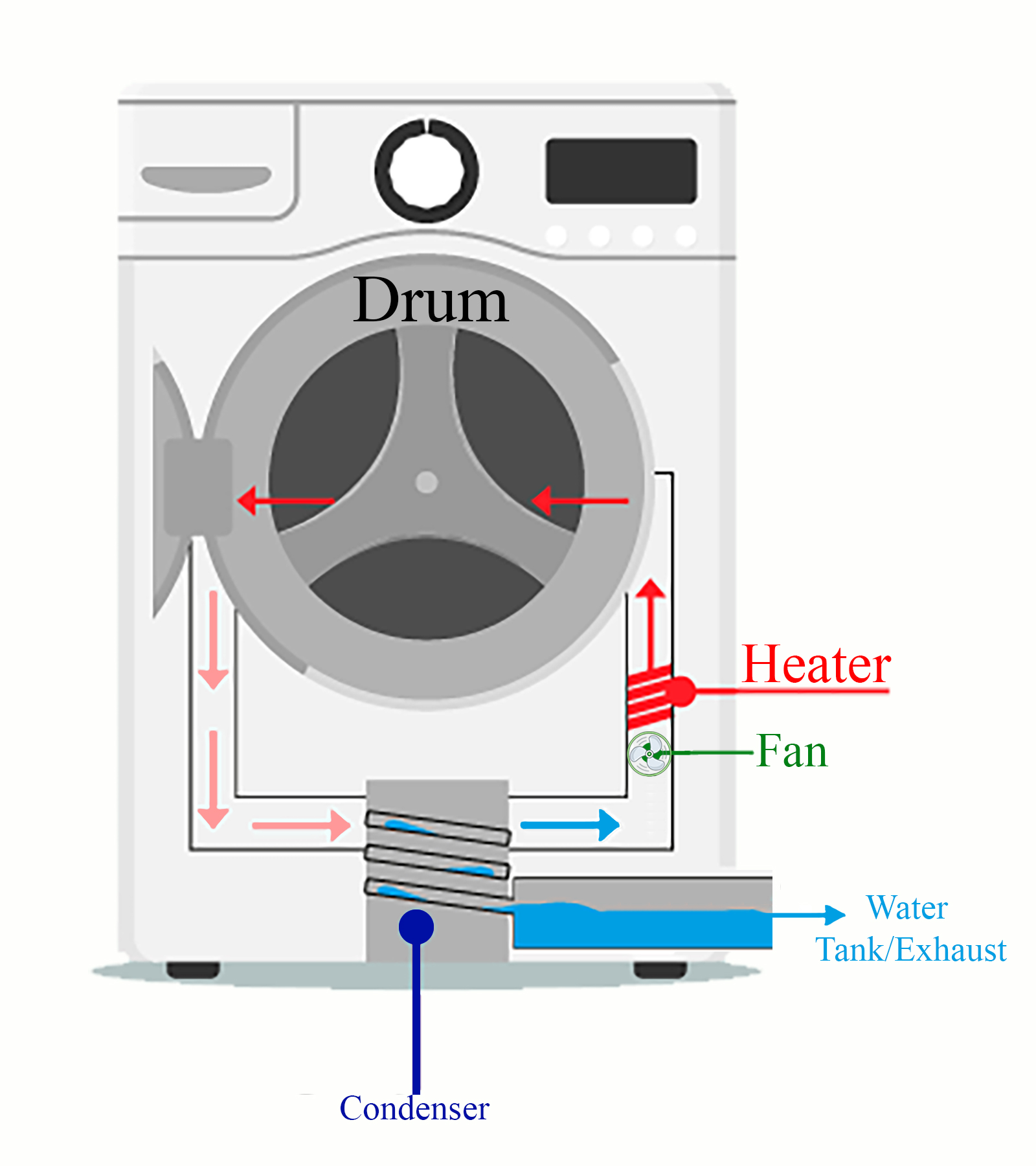}
\caption{Schematic of a drying cycle of the condenser dryers}
\label{CondenserDryer1}
\end{figure}

The clothes dryer works on the basis of flowing hot gas into the drum, where textiles are located, to remove the water content. Tumble dryers are divided into two types based on the cycle and energy source \citep{novak2019influence}. In the so-called \textit{vented dryers}, the drying takes place in an open cycle where hot air flows through the drum and exhausted the environment. Whilst, in the \textit{condenser dryers}, there is a close cycle of the gas in the machine which absorbs the water content from the textile and then loses humidity in the condenser. The schematic of these systems is shown in Fig. \ref{vented1} and \ref{CondenserDryer1}. 

\begin{figure}
\centering  \vspace*{0.5cm} \hspace*{0cm}
\includegraphics[scale=0.6]{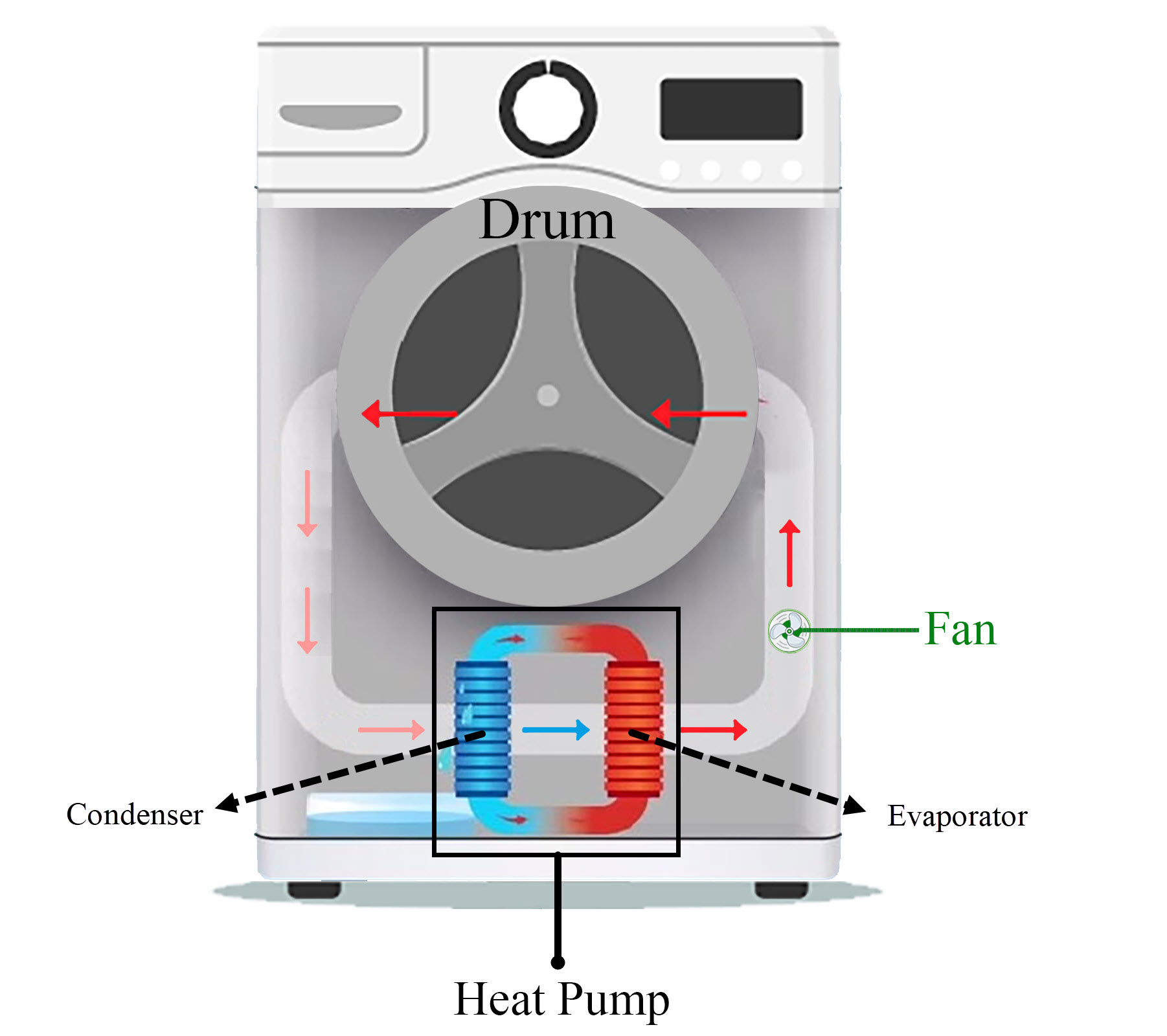}
\caption{Schematic of a drying cycle of the heat pump dryers}
\label{HeatpumpDryer}
\end{figure}

The vented dryers utilize hot gas or hot air provided by a gasifier or an electrical heater, respectively. Whereas in condenser dryers, the electrical heater is the source of the energy. In modern dryer machines, namely \textit{heat pump dryers}, a heat pump cycle is inserted into the dryer which helps the heat recovery throughout the system \citep{heatpumppatent1}. The heat pump tumble dryer is considered as a complex system due to the interdependency between the main components. The system briefly consists of two closed cycles, namely, the heat pump cycle, in which the working fluid is generally R134a \citep{rasti2012enhancement,senthilkumar2019influence,siddiqui2020recent}, and the drying cycle involving air-water mixture. The exchange of heat between cycles takes place at a condenser and an evaporator \citep{tegrotenhuis2017modeling}, as shown in Fig. \ref{HeatpumpDryer}.

Numerical models for the heat pump system development have been reported since the 2000s \citep{saensabai2003effects, sarkar2006transcriticala, sarkar2006transcriticalb, pal2008calculation, wang2010investigation, novak2011refrigerant}. Generally, a combination of experimental and numerical methods was utilized to investigate heat pump dryers, and more importantly, the portion of the numerical part increases through time \citep{svenssonflow, brandt2019exergetic}. 

The vented and heat pump dryers are installed in a separate machine from that of the washer machine, and hence, require more installation space. However, typical condenser dryers can be embedded in the washer machine, and therefore, are ideal for domestic applications. This specific machine in which both the washing and drying take place is called a \textit{washer-dryer}. In these machines, the gas mixture condenses via environment air or water as the coolant. Among all the mentioned drying systems, the condenser dryer is a promising method in which the drying needs fewer energy demands, drying time and no extra installation space. As a result, it accounts for a great portion of the USA market \citep{NorthAmericaHeatPumpDryers}.

Fig. \ref{condenserABC} shows the drying cycle in the condenser dryers. There are four main modules during the cycle; \textit{fan}, \textit{heater}, \textit{drum} and \textit{condenser}. Fan circulates the gas mixture through the cycle and has no effect on the thermodynamic properties of the mixture. However, the mixture undergoes severe thermodynamic changes in the other modules. Generally, the heaters are placed in the fan duct, and hence, it is conventional to treat them as a single module. Therefore, one can consider three separate sections after each module through the drying cycles. These sections are shown with \textit{a}, \textit{b} and \textit{c} symbols in Fig. \ref{condenserABC}. The main thermodynamic parameters related to the drying cycle are \textit{dry bulb temperature}, \textit{absolute humidity (humidity ratio)} and \textit{relative humidity}. Based on Fig. \ref{condenserABC}, in ideal conditions, the changes in thermodynamic parameters in each module are as follows:

\begin{itemize}
  \item Fan/Heater (\textit{a} $->$ \textit{b}): The dry bulb temperature increased. As no phase change occurs, the mass of the vapour and air, i.e. absolute humidity, is constant. The capability of the vapour absorption of the air is related to the temperature. Hence, the relative humidity decreases. As shown in Fig. \ref{Psycho0.jpg}, air moves through a horizontal line from point \textit{a} to \textit{b}.
  \item Drum (\textit{b} $->$ \textit{c}): As hot air flows through the wet clothes, it evaporates the water content and consequently cools down. Therefore temperature decreases. However, both absolute and relative humidity increase until they reach saturation conditions. This corresponds to the dotted path from \textit{b} to \textit{c} in Fig. \ref{Psycho0.jpg}.

Generally speaking, in the ideal conditions of a dryer, fully saturated air is expected at the outlet of the drum, i.e. relative humidity of 100\% \citep{lenaPhDThes, DrySchem}. However, in an actual machine, the outlet gas mixture may not reach saturation, especially at the middle or final stages of the drying period \citep{somdalen2018theoretical}. This is corresponding to the solid line from \textit{b} to \textit{c} in Fig. \ref{Psycho0.jpg}.

  \item Condenser (\textit{c} $->$ \textit{a}): Finally, in the condenser, the saturated gas cools down and its water content condenses. This happens at constant relative humidity conditions, shown with the path \textit{c} to \textit{a} in Fig. \ref{Psycho0.jpg}.

Finally, changes in the thermodynamic variables of the air in each module are listed in Table \ref{thermo1}.

 \begin{figure}
 \centering \vspace*{0.5cm} \hspace*{0cm}
 \includegraphics[scale=.5]{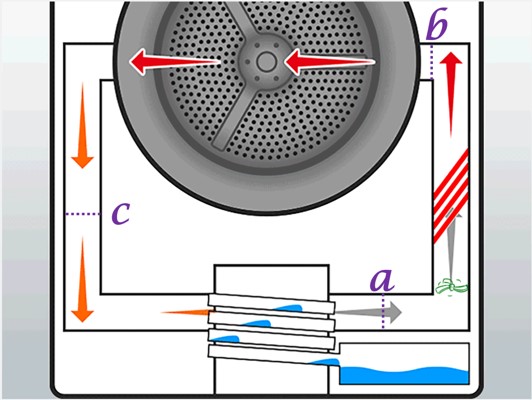}
 \caption{Different sections in the drying cycle of the condenser dryer.}
 \label{condenserABC}
 \end{figure}

 \begin{figure}
 \centering \vspace*{0.5cm} \hspace*{0cm}
 \includegraphics[scale=.25]{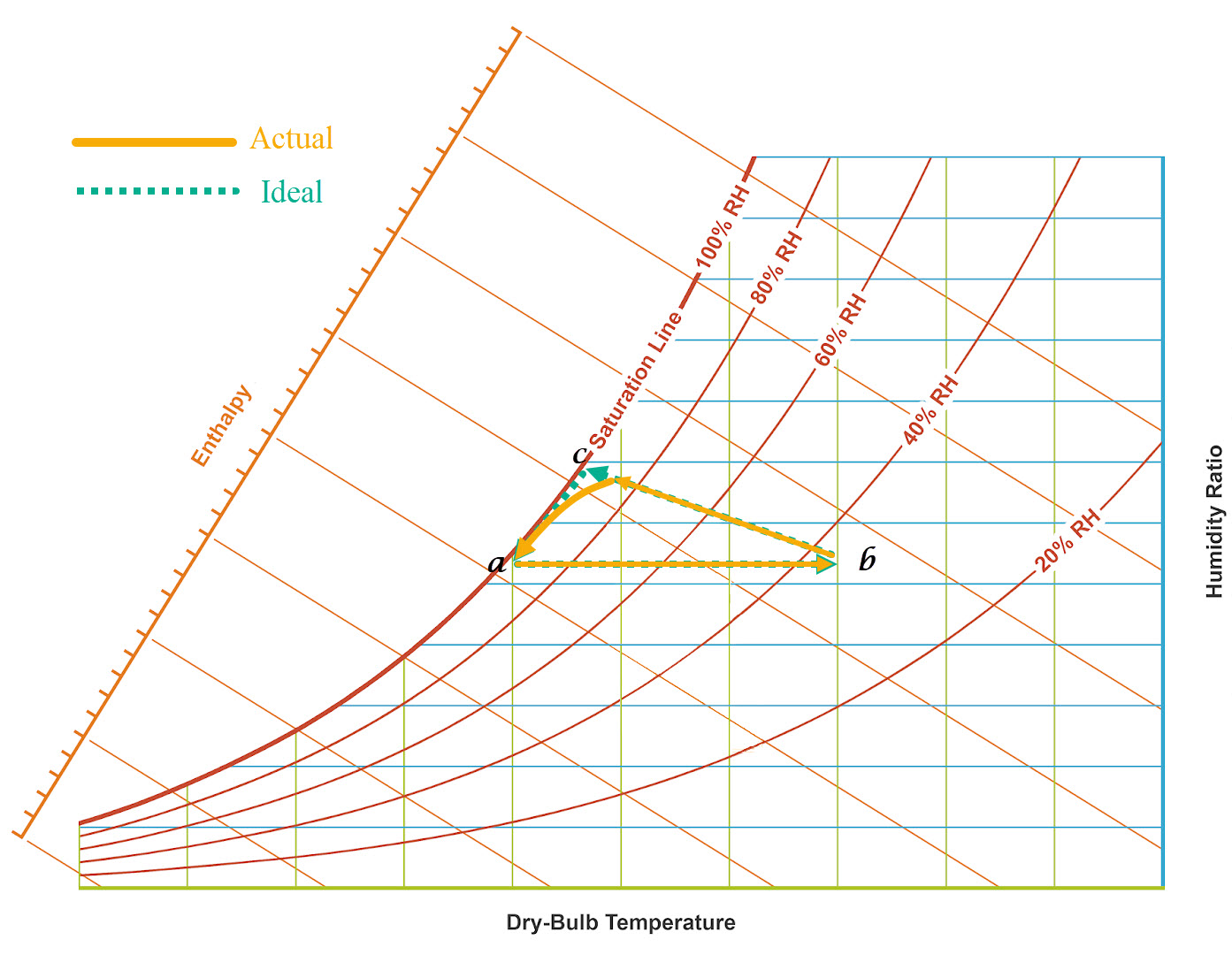}
 \caption{Ideal drying cycle in the psychometric chart.}
 \label{Psycho0.jpg}
 \end{figure}  
 
\begin{table*}
\centering
\caption{Thermodynamic changes of the mixture in the cycle}
\label{thermo1}
\renewcommand{\arraystretch}{2}
\begin{tabular}{|c|c|c|c|}
\hline
\textbf{Module} 			& \textbf{Fan/Heater} 			& \textbf{Drum} 			& \textbf{Condenser} 		\\ \hline
\textbf{Path} 			& \textbf{a $->$ b} 			& \textbf{b $->$ c} 		& \textbf{c $->$ a} 		\\ \hline
\textbf{Dry Bulb Temp. }         & $T_{dry}$$\uparrow$               & $T_{dry}$$\downarrow$       & $T_{dry}$ $\downarrow$     \\ \hline 
\textbf{Relative Humidity}            & RH$\downarrow$               & RH$\uparrow$                  & RH$\simeq$                \\ \hline
\textbf{Absolute Humidity}            & AH$\simeq$                   & AH$\uparrow$                  & AH$\downarrow$            \\ \hline\end{tabular}
\end{table*}

\end{itemize}

\section{Physical Phenomena}
The analysis of the drying phenomenon is more complicated than that of heat or isothermal mass transfer alone \citep{choi2002drying,park2003drying}. During the drying of wet materials, heat and mass transfer occurs simultaneously in both the solids and in the boundary layer of the drying agent \citep{lee2006drying}. 

Specifically speaking, in a condenser dryer, the gas mixture and clothes both undergo several physical phenomena. Generally, while the drum and clothes in it are spinning with the aid of the motor, hot and almost dry air enters the drum. In the drum, hot gas mixes with the wet clothes and passes through the pores of textiles. During this process, all types of heat transfer and mass transfer mechanisms occur. There exist conductive and convective heat transfers between clothes and clothes to the air, respectively, while radiation takes place between components all over the dryer. The mass transfer contributes to the evaporation of clothes' water content into the air and condensation of the air mixture, respectively, through the drum and condenser. Adding to mentioned complexities are the transient effects introduced by the on/off heater to maintain constant drying temperature along with the random motion of clothes in the drum. The phenomena involved in the drying process and the module in which they take place are listed in Table \ref{tab-phyphen}.

\begin{figure*} 
\centering  \vspace*{0cm} \hspace*{0cm}
\includegraphics[scale=0.35]{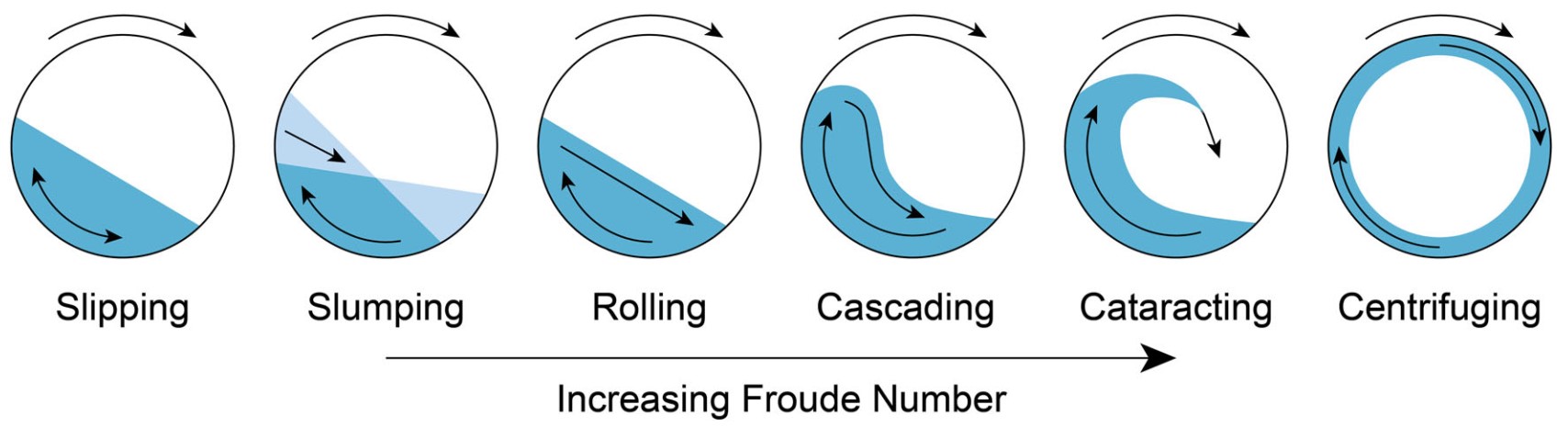}
\caption{Flow regimes in rotating drums with increasing Froude number \citep{jones2022dynamics}}
\label{ImagingTechFroude}
\end{figure*}

Moreover, Jones et al. identified different two-phase regimes in the drum during the drying cycle \citep{jones2022dynamics}. Accordingly, there exist six flow regimes as a function of Froude number in the drum, as shown in Fig. \ref{ImagingTechFroude}, which demonstrate the range of movements experienced by items during tumbling. It was concluded that, as the rotation rate increases, more chaotic flow and turbulence are expected resulting in a fully annular flow at high angular velocities. This research did not manage to relate those regimes to the drying properties.

\begin{table*}[]
\centering
\caption{Physical phenomena}
\label{tab-phyphen}
\renewcommand{\arraystretch}{2}
\begin{tabular}{|c|c|}

\hline
\textbf{Phenomenon}                          & \textbf{Module}             \\ \hline
Multi-component flow (air and vapor)          & Fan/Heater, Drum, Condenser \\ \hline
Multi-body dynamics (MBD) (clothes rotation) & Drum                        \\ \hline
Heat transfer (convection and radiation)     & Fan/Heater, Drum, Condenser \\ \hline
Heat Transfer (conduction between clothes)   & Drum                        \\ \hline
Variable thermophysical properties          & Fan/Heater, Drum            \\ \hline
Multiphase flow (gas and water)              & Drum, Condenser             \\ \hline
Mass Transfer (evaporation and condensation) & Drum, Condenser             \\ \hline
\end{tabular}
\end{table*}

\section{Methodologies}

\begin{figure} 
\centering  \vspace*{0cm} \hspace*{0cm}
\includegraphics[scale=0.25]{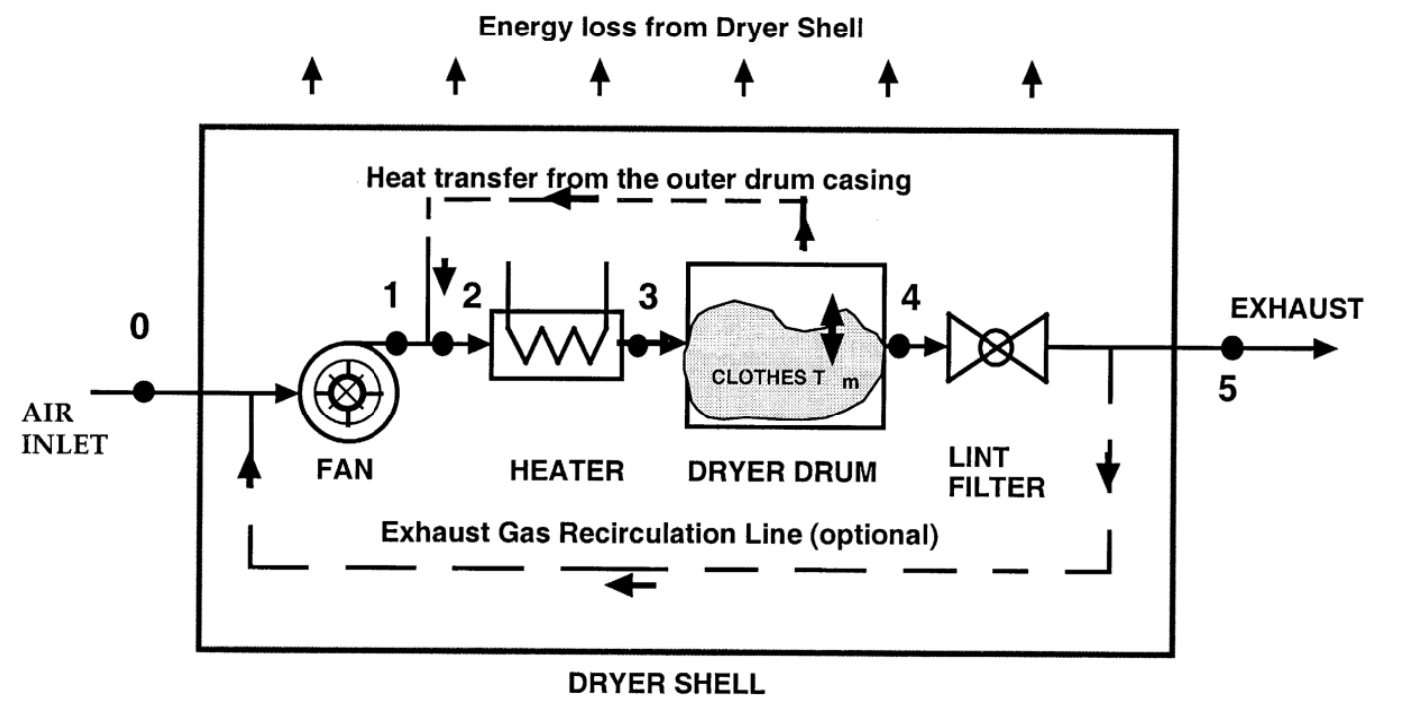}
\caption{Box of the control volume around a semi open-loop vented dryer \citep{deans2001modelling}}
\label{deansnodes}
\end{figure}

The involving modules in a tumble dryer are common equipment in various industries, and thus, there exists plenty of literature dealing with a fan/heater, rotating drum or condenser \citep{shen2013performance, yousef2014experimental, enteria2015performance, mensah2017energy, zhang2019numerical, khalsa2021experimental, zolotarevskiy2022modelling}. Moreover, as the drum is the most important module in the cycle in which the drying takes place, if possible, examining this module alone can provide us with the necessary parameters for evaluating drying efficiency. But unfortunately, derivation of the drying cycle must be carried out by considering involving components altogether, and hence, the so-called modular approach is not applicable to tumble dryers. This is due to the fact that the operating condition of each module may have an effect on the working condition of other modules. For instance, the gas flow rate in the fan channel corresponds to the pressure loss through all other modules. A similar explanation is true for mass transfer rates in the drum or condenser. As a consequence, even for determining the performance of a single component, one needs to consider all the modules. Therefore, the literature has preferred to consider a control volume approach encompassing all the components and solving the mass and energy conservation for the whole system, such as the one shown in Fig. \ref{deansnodes}, proposed by Deans for the open-loop vented dryer, which will be discussed later on.

There were a little number of publications in the open literature regarding laundry drying technologies before the 90s \citep{conde1997energy}. In addition, earlier studies were directed more at the process of drying textiles during the manufacturing stage \citep{conde1997energy}. Primary research on the drum module in the tumble dryers was started with observations and the use of experimental measurement \citep{fontana2016mathematical,akyol2015simulation}. Hekmat and Fisk published one of the earliest studies in which practical methods including reduction of the airflow rate and recirculation of the exhaust air, to improve the efficiency of the drying \citep{hekmat1984improving}, which can result in 33\% energy saving at the cost of severe drying time increase \citep{boudreaux2021measurement}. Spruit measured temperature and humidity change in the drum to qualify drying performance \citep{spruit1988model}. They found out that an increase of the textile load of 20\% above the nominal resulted in an increase in the energetic drying efficiency of the tumbles by 66\%. They also observed the effect of variation of the inlet variables on the drying performance. One of the first attempts in optimizing the condensation-type clothes dryer was accomplished by Kadoya et al. \citep{kadoya1984condensation}. They tried to optimise the fan of the heat exchanger which resulted in a simplified system with less weight and consequently cheaper equipment. This work was further developed by \citep{kim1999optimum} to optimize fan performance and decrease the noise with the aid of the Taguchi design of experiment (DOE) method \citep{roy2001design}.

The main issue with early-stage research was the lack of theoretical focus. They follow a rather similar routine in measuring the important parameters such as temperature and humidity indices in the drying period and observe the sensitivity of those variables to inlet conditions. They even succeeded in optimizing and predicting drying workflow. For instance, Lemaire et al. reported the effects of the systematic variation of the drying process parameters (air inlet temperature, load ratio, recirculation rate, ventilator capacity and drum rotation speed) and established an empirical model of the drying process for their gas-heated tumble \citep{conde1997energy}. However, with the lack of an appropriate mathematical model, it is hard to generalize those findings to other machines.

\subsection{Lambert Model} 
It is hard to understand the relation between machine parameters without inspecting theoretical aspects of physical phenomena in the drying cycle. As a consequence, there is no proper interrelation between studies before the 90s which were merely experimental. Mathematical modelling of the drying cycle for fabrics requires prior knowledge of background physics which enables process design and minimization of energy costs \citep{kahveci2007transport}.

The very first study on drying machines with a focus on mathematical aspects was done by \cite{lambert1991modeling}. Their goal was to derive a heat and mass transfer model of the drying process in a textile tumble dryer based on Chilton–Colburn analogy (also known as the modified Reynolds analogy \citep{ccengelheat}) which relates heat and mass transfer phenomena \citep{ashrae1985ashrae}. It should be mentioned that this approach has been widely used in the literature of drying for porous media, semi-conductors, food, agriculture, Condebelt and distillation \citep{tsilingiris2012combined,tsilingiris2015theoretical, collazo2018mathematical,gu2018effects,brummans2019modeling, he2020modeling,buysse2021dryer,sabau2020evaporation, nienke2021experimental,chasiotis2021evaluation}. In this regard, they estimated the mass evaporation rate in the drum with the following relationship:

\begin{equation}
\dot{m}_{\mathrm{ev}}=k A \left(P_{\mathrm{t}}-P_{\mathrm{x}}\right) = k A (\omega_{s}-\omega_{out}) \;, \\
\label{equlam}
\end{equation}

where $k$ and $A$ are the evaporation coefficient and surface area of the clothes, respectively. Conventionally, the multiplication, $k A$, which is known as the area-mass transfer coefficient, is considered as a whole. All the physical processes in the drum can be described by the area-mass transfer coefficients which will be described in more detail further on. 

$P_{\mathrm{t}}$ and $P_{\mathrm{x}}$ are partial vapour pressure at the textile's surface and the drying air (bulk air), respectively. $P_{\mathrm{t}}$ is estimated as follows:

\begin{equation}
P_{\mathrm{t}}\left(x_{\mathrm{t}}, T_{\mathrm{t}}\right)=P_{\mathrm{sat}}\left(T_{\mathrm{t}}\right) a\left(x_{\mathrm{t}}\right) \;,\\
\label{psat}
\end{equation}

in which $P_{\mathrm{sat}}$ is a function of the temperature. The water activity, $a$, is a function of the moisture content of the medium to be dried. Although the water activity also depends in a complicated way on the temperature and direction of the process (e.g. drying rather than wetting), the application of approximate curves is satisfactory on the level of aggregation. Describing the curve of the $a$ is commonly known as sorption-isotherms in the literature. 

\begin{figure} 
\centering  \vspace*{0cm} \hspace*{0cm}
\includegraphics[scale=0.75 ]{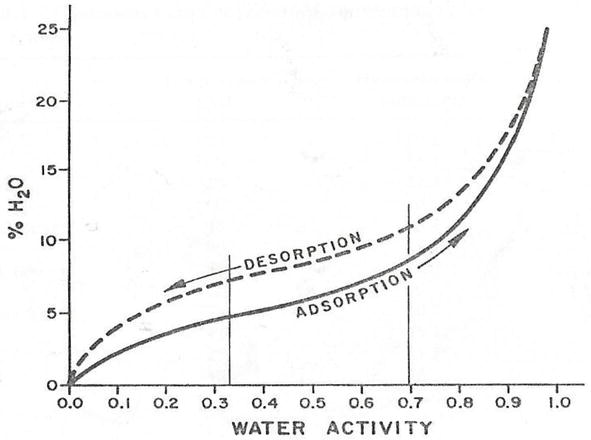}
\caption{ A typical moisture sorption isotherm curve \citep{aviara2020moisture}}
\label{soriso}
\end{figure}

The sorption-isotherm curve demonstrates the relationship between the textile's water content and the equilibrium relative humidity (water activity). Put differently, for each humidity value, a sorption-isotherm indicates the corresponding water content value at a given constant temperature \citep{glaskova2009moisture}. An increase in $a$ is usually accompanied by an increase in water content but in a non-linear fashion \citep{ghimire2017basic}. These curves are merely determined with experiments and cannot be explicitly derived by calculations. A typical sorption-isotherm curve is shown in Fig. \ref{soriso}. As depicted, the behaviour of the material in evaporation and condensation does not follow the same line. Generally, the behaviour of this curve can be of any type which varies based on the fabric's material. Another typical behaviour in sorption-isotherm curves is the occurrence of hysteresis \citep{jiang2019modelling}. 

\begin{figure} 
\centering  \vspace*{0cm} \hspace*{0cm}
\includegraphics[scale=0.4 ]{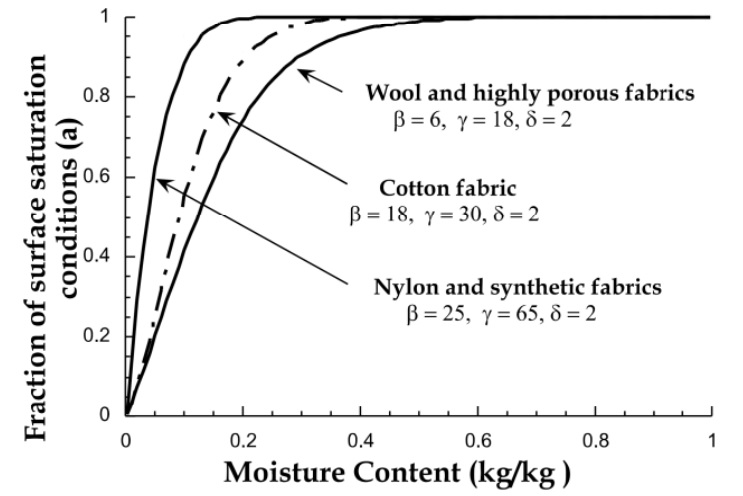}
\caption{ Desorption isotherms for different fabrics \citep{aviara2020moisture}}
\label{soriso2}
\end{figure}

For the drying process, one only needs to account for the desorption process (evaporation). The first attempt to describe the desorption-isotherm of various textiles is done by \citep{Krischer1978} which is regarded as the primary reference in this field. Lambert regressed the following relationship for the estimation of activity coefficient \citep{lambert1991modeling}:

\begin{equation}
a=1-\frac{\beta^{*} M_{\mathrm{c}}+\delta}{1+\delta^{\gamma * M_{\mathrm{c}}}}
\label{waterActivity}
\end{equation}

where $ M_{\mathrm{c}}$ is the moisture content of clothes and other coefficients depend on the fabrics' type. Fig. \ref{soriso2} presents coefficients of activity for various fabrics, i.e. a highly porous fabric such as wool, and low porosity fabrics such as nylon and cotton. Based on the figure, the surface of the cotton fabric will be saturated during most of the drying cycle. It is notable that it is only in the later stages of the drying process that the desorption isotherm of the different fabrics will become significant.

It should be noted that through time, more complex relations were proposed for the water activity for the fabric drying process. For instance, Yi model which relies on the water content and clothes temperature is as follows \citep{yi2015new}:

\begin{equation}
a\left(x, T_{c l}\right)=1-e^{\left(-\alpha . x(t) . h_{f g} / k_{B} . T_{c l}\right)} \; .
\end{equation}

\begin{figure} 
\centering  \vspace*{0cm} \hspace*{0cm}
\includegraphics[scale=0.55 ]{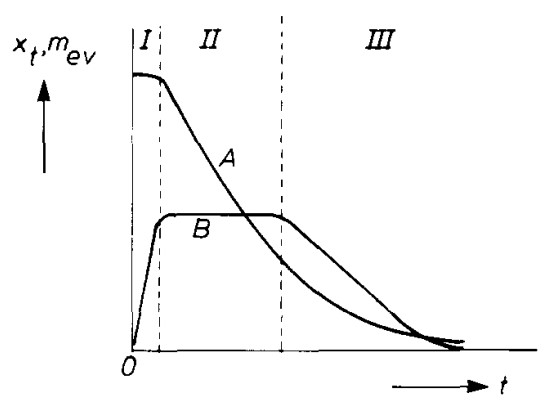}
 \caption{Schematic of the drying process: transient (I), constant rate (II) and falling rate (III) period. (A) Moisture contents, $x_{t}$, versus time, $t$; (B) rate of evaporation, $\dot{m}_{ev}$, versus time \citep{lambert1991modeling}}
\label{evapRate1}
\end{figure}

With the aid of Eq. \ref{equlam}, Lambert reported the evaporation rate during a dryer cycle, with which he could classify it into three stages. In general, the drying process starts with a transient change from the initial to a stationary state. Although some local evaporation may occur in the drum, a huge portion of the heater energy is dedicated to increasing the clothes' temperature. This stage, also known as the \textit{transient}, is depicted in Fig. \ref{evapRate1}. In this period, the moisture content is almost constant and the evaporation rate starts to increase. 

In the next stage, i.e. \textit{constant rate} or \textit{stationary state}, the water content decreases monotonically. In other words, the input power is primarily consumed to evaporate the clothes' water content. It should be noted that this stage happens if the moisture content of the textile is sufficiently high to keep the water activity approximately equal to unity for at least some period of time, which is a common condition in dryer machines. Furthermore, all process parameters will tend to remain constant, at this stage, especially the rate of evaporation. The condition of unity for the water activity means that the evaporation rate does not depend on the desorption-isotherm. 

When the moisture content, $x_{t}$, becomes sufficiently small, the water activity deviates from unity. This occurs because the assumption of a free water surface is no longer valid. Decreasing the water activity is accompanied by a decrease in the rate of evaporation. That in turn leads to an excess power that is directly consumed to heat the textile and drum. Thus, the reduction in water activity will be compensated partially by an increase in $P_{sat}$ (Eq. \ref{psat}). This leads to the \textit{falling-rate} stage, which is characterised by an increasing temperature of the exhaust air and a decreasing rate of evaporation, $\dot{m}_{ev}$  \citep{lambert1991modeling}. 

Back to Eq. \ref{equlam}, the remaining unknown is the area-mass transfer coefficient, $kA$. Still, there is no universal mathematical model to address this parameter which corresponds to the complexity of the mass transfer mechanism. Hence, this coefficient has been treated differently in the literature. Some authors utilized the concept of analogy between heat and mass transfer to address this variable \citep{huang2019simulation}. Modified Lewis number can relate convective heat transfer and mass transfer coefficients which is defined as follows \citep{MollyPhDThes}: 

\begin{equation}
\mathrm{Le^{\prime}}=\frac{h}{\lambda C_{P}}
\label{deansEq5}
\end{equation}

The approximate ranges of convective heat transfer coefficient, $h$, are listed in Table \ref{tab:h}. Although at the first sight, a coefficient based on the phase change should be considered, by taking a closer look, the heat transfer actually occurs between water in the clothes and the saturated air layer above. Thus, it must be treated as forced convection, and therefore, the value of $h=20\sim100$ seems more applicable. However, the appropriate value is not known yet. Considering the high difference between the top and bottom range of $h$, a non-precise value can dramatically change the evaporation rate. It is worth mentioning that the determination of effective clothes area, $A$ in Eq. \ref{equlam}, has its complexities and is treated differently through the literature \citep{huang2019simulation, adapa2002performancea, adapa2002performanceb}.

\begin{table}[]
\centering
\caption{Conductive heat transfer coefficient \citep{bergman2011introduction}}
\label{tab:h}
\begin{tabular}{|c|c|}
\hline
\textbf{Process}                         & \textbf{h(W/m2·K)} \\ \hline
Natural Convection                       &                    \\
Air                                      & 1$\sim$10               \\
Water                                    & 200$\sim$1000           \\ \hline
Forced Convection                        &                    \\
Gas                                      & 20$\sim$100             \\
High pressure vapor                      & 500$\sim$35000          \\
Water                                    & 1000$\sim$1500          \\ \hline
Phase Change of Water &                    \\
Boiling                                  & 2500$\sim$3500          \\
Vapor condensation                       & 5000$\sim$25000         \\ \hline
\end{tabular}
\end{table}

The mentioned complications convinced authors to utilize a constant value of area-transfer coefficient based on the experimental data \citep{deans2001modelling}. For instance, \citep{bengtsson2014performance} considered a value of $0.036 kg/s Bar$. Although this approach provides more degrees of precision as it is derived experimentally, considering a constant value for the coefficient is still unfavourable. As mentioned earlier, especially at the final stages of drying, the evaporation rate declines. Hence, this method is not valid if the whole drying period is being simulated. Hence, literature tended to employ this assumption in the simulation of heat pump tumble dryers. Basically, due to the complexity imposed by considering the heat pump cycle, the utilization of a constant area-transfer coefficient is convincing. 

\begin{figure} 
\centering  \vspace*{0cm} \hspace*{0cm}
\includegraphics[scale=0.28 ]{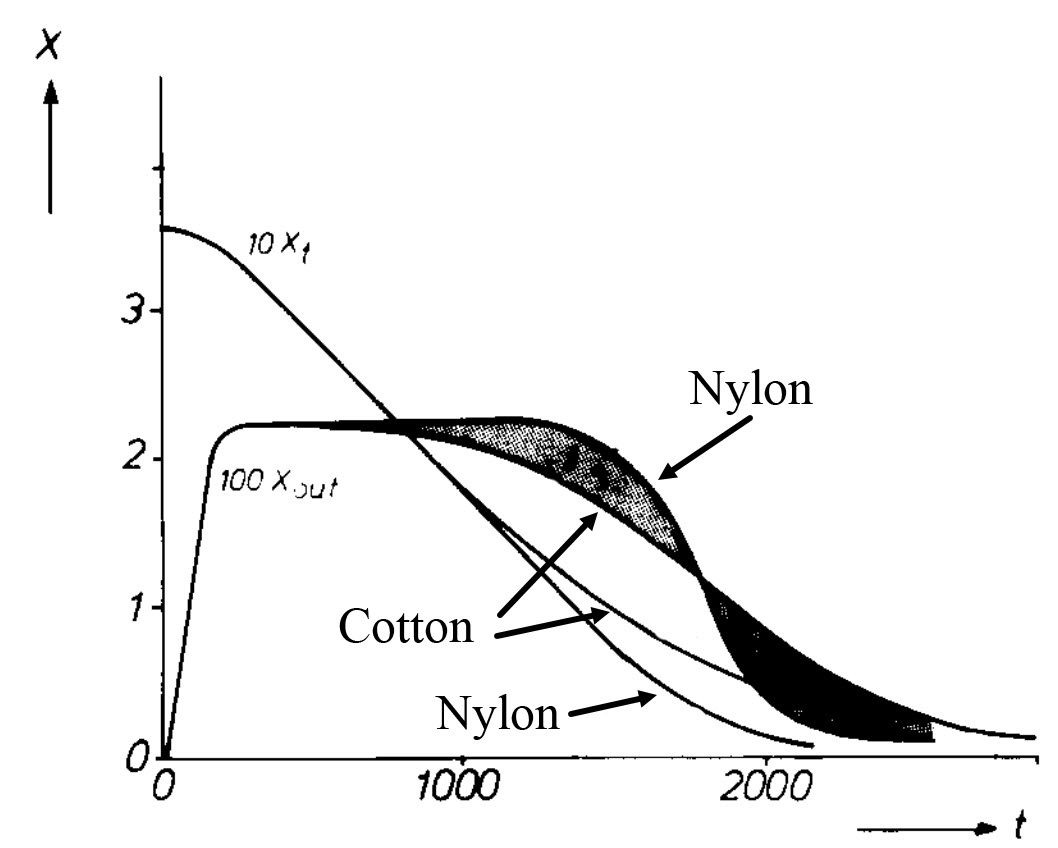}
 \caption{Schematic of the drying process for two types of fabrics \citep{lambert1991modeling}}
\label{evapRate2}
\end{figure}

Back to the Lambert model, Eq. \ref{equlam}, he was able to calculate the evaporation rate of nylon and cotton fabrics over time. The resulting evaporation curve along with moisture content is depicted in Fig. \ref{evapRate2}. The moisture content of the two fabrics has a similar behaviour during \textit{transient} and \textit{constant rate} stages. There exists a deviation at the beginning of the \textit{filling rate} stage where the water activities are different and far from unity.

There exist other methodologies for deriving a relationship for evaporation in tumble dryers. However, the drying models can be challenging and computationally expensive \citep{turner1996mathematical}. Ng and Deng developed a semi-empirical model for clothes drying based on a termination control method \citep{ng2008new}. However, their model covers only the last constant-rate and falling-rate periods of the drying cycle. Haghi employed a mathematical model with the rate equation for mass transfer between the moisture in the fabric fibres and the moisture in the air within the fabric pores \citep{haghi2001mathematical, nordon1967coupled}. These types of models provide some local insights but are too complex and computationally expensive to be used for optimization purposes. Moreover, they require an immense amount of input regarding thermophysical properties and heat and mass transfer coefficients, and therefore, did not gain much attention.

\subsection{Deans Model}
Due to the simplicity and feasibility of the Lambert model, several studies have utilized or developed the model to predict the drying cycle \citep{ bejan2004porous, lee2006dryingHandbook, stawreberg2013potential}. Deans developed a mathematical model to analyze the drying performance of the dryer system based on the balances of mass and energy \citep{deans2001modelling, jian2017drying}. The focus of the model was the heat and mass transfer flows within the dryer drum by utilizing total mass and energy balances for the air/vapour mixture at nodes throughout the dryer. The location of these nodes along with a description of the air/vapour flow and energy flow routes is shown in Fig. \ref{deansnodes}. This restructured model enables him to track energy inputs and losses from the dryer along with the implementation of the airflow recirculation. The use of total energy balances for the air/vapour mixtures permits the identification of the components which are unimportant and can therefore be ignored. The modular construction used in the model enables variation to the machine inputs such as heater power to be made during the cycle. This approach is important when potential energy-saving measures are investigated \citep{deans2001modelling}. 

Frankly, the main focus of Dean's model is the dryer drum since the efficient operation of this component determines the mass and energy flows throughout the complete unit. The main simplification used in the model is the assumption of identical temperature between fabrics, the water content of the clothes and the drum walls. As the most difficult variable to measure in the drying cycle is the clothes temperature \citep{gluesenkamp2018experimental}, they expressed this variable based on the energy balance in the drum as follows \citep{deans2001modelling}:

\begin{equation}
\frac{\mathrm{d} T_{\text {clo}}}{\mathrm{d} t}=\frac{H_{\text {in}}-H_{\text {out}}-H_{\text {evap}}-\dot{Q}_{\text {loss}}}{\left(\dot{m} C_{p}\right)_{\text {dry clo}}+\left(\dot{m} C_{p}\right)_{\text {moist}}+\left(\dot{m} C_{p}\right)_{\text {d}}} \;,
\label{deanEq1}
\end{equation}

which can be rewritten as follows \citep{deans2001modelling}:

\begin{equation}
\frac{T_{\text {out}}-T_{\text {clo }}-\left(\frac{H_{\text {evap}}-\dot{Q}_{\text {loss }}}{h A}\right)}{T_{\text {in}}-T_{\text {clo}}-\left(\frac{H_{\text {evap}}-\dot{Q}_{\text {loss }}}{h A}\right)}=\mathrm{e}^{-\left[\frac{h A}{\left(m_{\mathrm{a}} C_{p}\right)_{\text {average}}}\right]}
\label{deanEq2}  \;,
\end{equation}

which equal to \citep{lee2019rationally}:

\begin{equation}
T_{\text {c}+d t} = T_{\text {c}}+\\
\left[\frac{\dot{m}_{a}\left(h_{\text {in }}-h_{\text {out }}\right)-\dot{Q}_{\text {loss }}+\dot{m}_{\text {evap }} h_{\text {evap }}}{\dot{m}_{\text {dry clo}} C p_{\text {dry clo}}+\dot{m}_{\text {dry clo}} C p_{moist}+\dot{m}_{d} C p_{d}}\right] \Delta t  \;.
\label{deanEq22}
\end{equation}

\vspace{0.2cm}

Eg. \ref{deanEq22} is a well-known equation for estimating the temperature evolution in the drum which has been implemented in several studies \citep{ bejan2004porous, lee2006dryingHandbook, stawreberg2013potential}. Specific heat capacities are dependent on the clothes and machine materials and are generally assumed to be constant as they have insignificant changes in terms of temperature. $h_{evap}$ is the enthalpy of vaporization (the latent heat of evaporation) and is a function of temperature. $\dot{m}_{a}$ is the dry air mass flow rate entering the drum, which is determined experimentally. 

$Q_{\text {loss }}$ is the convective heat transfer loss from the outside walls of the tumble dryer to the ambient air, which is a function of machine dimensions and temperatures of the gas mixture and ambient. This term varies in each step of the cycle as the temperature of the gas mixture changes. However, the highest heat loss occurs in the heater module in which the gas mixture has the highest temperature. For the sake of simplicity, Deans took it as a function of the temperature in the heater module. Finally, the heat loss term takes the following form \citep{bengtsson2014performance}:

\begin{equation}
{Q}_{\text {loss }}=\mathrm{Nu} \cdot \frac{\lambda_{\text {a }}}{L_{th}} \cdot A_{ext}\left(T_{hater}-T_{amb}\right) \;,\\
\label{QLoss}
\end{equation}

where $\lambda_{\text {air }}$ is the air thermal conductivity. ${L_{th}}$ and $A_{ext}$ are the shell thickness and external surface area of the dryer, respectively. Nusselt number is the ratio of convective to conductive heat transfer mechanisms and is defined as \citep{bergman2011introduction,granryd2005secondary}:

\begin{equation}
\mathrm{Nu}=0.13 \cdot(\mathrm{Gr} \cdot  {Pr})^{\frac{1}{3}} \quad 10^8<\mathrm{Gr} \cdot  {Pr}<10^{12}
\label{NuEq}
\end{equation}

$\mathrm{Gr}$ and $  {Pr}$ are Grashof and Prandtl numbers, respectively, and are defined as follows:

\begin{equation}
\mathrm{Pr}=\frac{\nu}{\alpha}=\frac{\text { Momentum diffusivity }}{\text { Thermal diffusivity }} \\
\label{PrEq}
\end{equation}

\begin{equation}
\mathrm{Gr}=\frac{g \beta\left(T_{heater}-T_{ambient}\right) {L_{th}}^3}{\nu^2}=\frac{\text { Buoyancy Force }}{\text { Viscous Force }} \\
\label{GrEq}
\end{equation}

where $g$ and $\beta$ are gravitational acceleration and the thermal expansion coefficient, respectively.

Deans proposed the following relation instead of that of Lambert (Eq. \ref{equlam}) for the evaporation rate:

\begin{equation}
\dot{m}_{\mathrm{evap}}=k A\left[a \omega_{\mathrm{M}}-\frac{\left[\omega_{\text {inlet }}+\omega_{\text {exhaust }}\right]}{2}\right] \;,
\label{deansEq3}
\end{equation}

and therefore:

\begin{equation}
\dot{m}_{\mathrm{a}}\left[\omega_{\text {inlet }}-\omega_{\text {exhaust }}\right]=k A\left[a \omega_{\mathrm{M}}-\frac{\left[\omega_{\text {inlet }}+\omega_{\text {exhaust }}\right]}{2}\right]
\label{deansEq4}
\end{equation}

For the sake of simplicity, some authors considered the humidity ratio of clothes equal to that of the inlet \citep{lee2019rationally}:

\begin{equation}
\dot{m}_{\mathrm{a}}\left[\omega_{\text {inlet}}-\omega_{\text {exhaust }}\right]=k A\left[a \omega_{\text{inlet}}-\frac{\left[\omega_{\text {inlet }}+\omega_{\text {exhaust }}\right]}{2}\right]  \;,
\label{deansEq4_2}
\end{equation}

where $\omega_{\mathrm{M}}$, $\omega_{\text {inlet}}$ and $\omega_{\text {exhaust }}$ are the humidity ratios of the saturated air at clothes surface and air at the drum inlet and outlet. $\dot{m}_{\mathrm{a}}$ is the dry air mass flow rate at the inlet. The only unknown in Eq. \ref{deansEq4} is the specific humidity at the outlet, $\omega_{exhaust}$, which can be determined through time and consequently, evaporation rate, $\dot{m}_{\mathrm{evap}}$ is calculated.

\begin{figure*} 
\centering  \vspace*{0cm} \hspace*{0cm}
\includegraphics[scale=0.5]{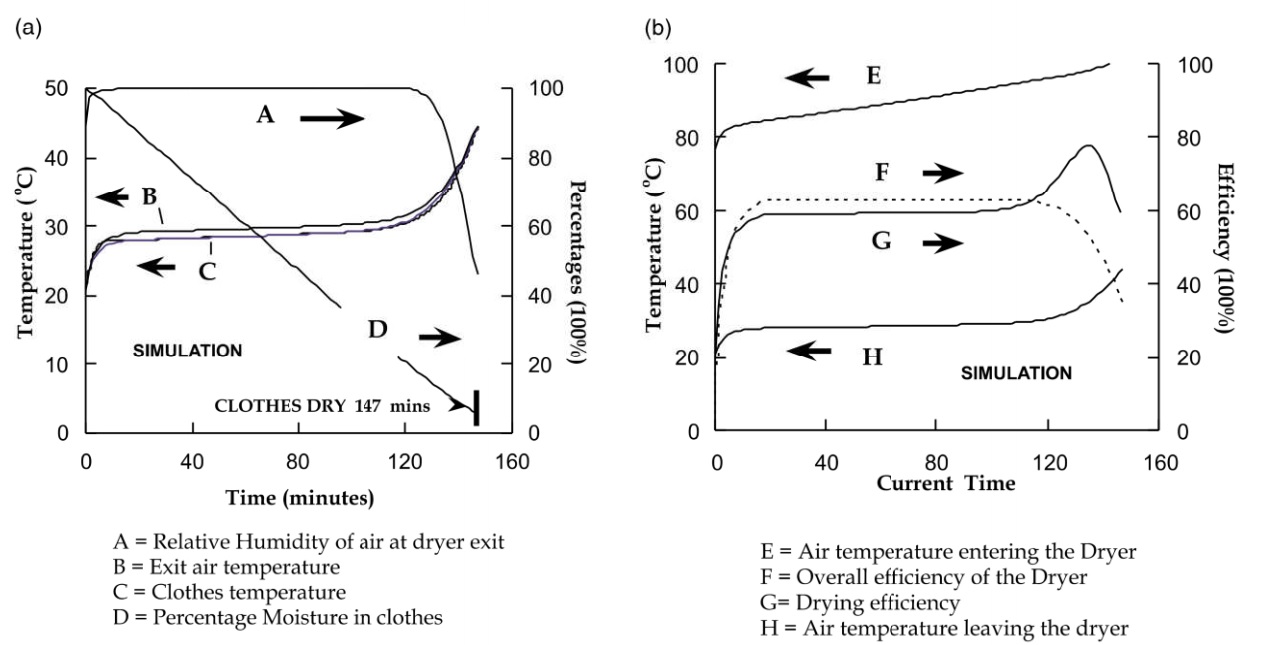}
 \caption{Results of numerical model of Deans \citep{deans2001modelling}}
\label{deansoutputs1}
\end{figure*}

By solving energy and mass balance equations, Deans predicted the evolution of the variables in the drying machine through time which is depicted in Fig. \ref{deansoutputs1}. Drum-specific parameters are shown on the left. The relative humidity (chart A) in the drum is predicted to be 100\% throughout almost the drying cycle except for the start and final stages (transient and falling rate stages, respectively). Whilst, the clothes and exit temperatures (charts B and C) show a similar alteration. The change of the moisture content (chard D) is comparable to what was described earlier (Fig. \ref{evapRate1}). On the right, other indices regarding the whole control volume (Fig. \ref{deansnodes}) are presented. 

\begin{figure*} 
\centering  \vspace*{0cm} \hspace*{0cm}
\includegraphics[scale=0.4]{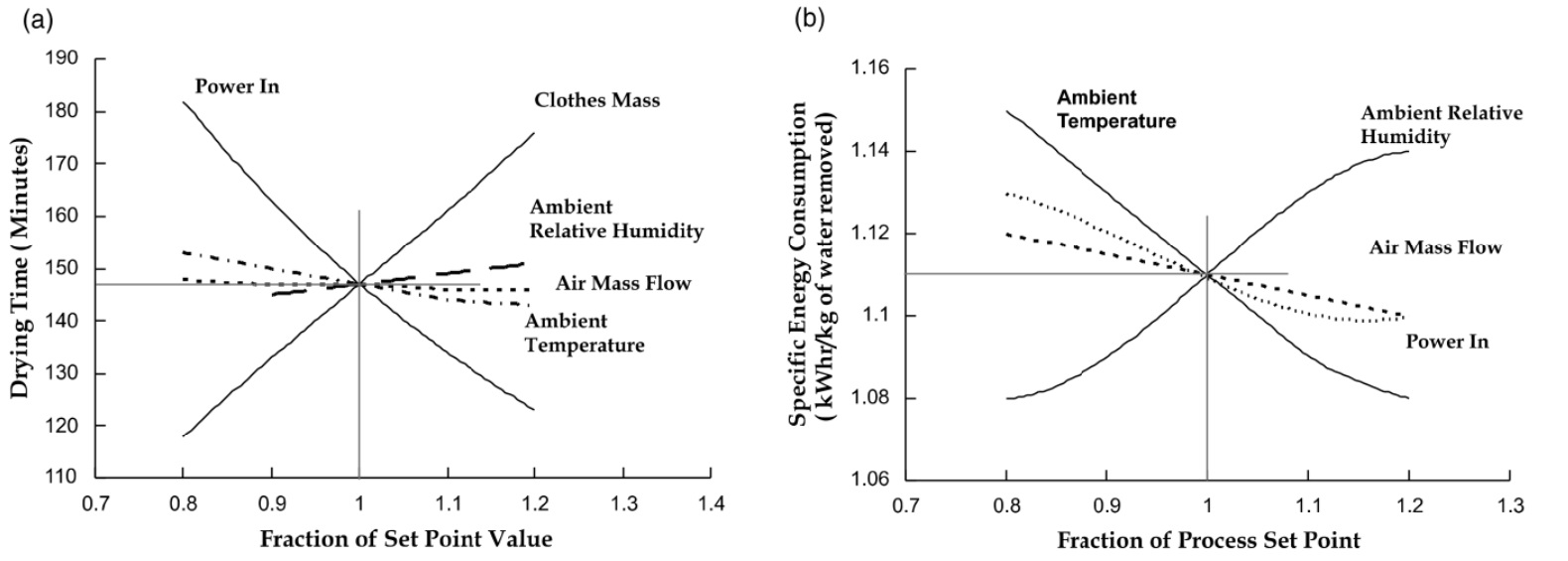}
 \caption{The sensitivity of drying times and process efficiency with design parameters \citep{deans2001modelling}}
\label{deansoutputs2}
\end{figure*}

Additionally, he has successfully carried out a sensitivity analysis on the drying time and efficiency numerically which is depicted in Fig. \ref{deansoutputs2}. Accordingly, by considering the current working status of the system as the set point (intersection point in the diagram), changes in the drying time and energy consumption in terms of the input power, clothes mass, ambient temperature and relative humidity, as well as inlet flow rate, were estimated numerically. Ambient humidity has a negative effect on drying time and energy consumption. While the increase in the ambient temperature, air flow rate and input power results in faster drying with less energy consumption. It should be noted that input power accounts for the fan, rotating drum and heater powers. 

\begin{figure}
\centering  \vspace*{0cm} \hspace*{0cm}
\includegraphics[scale=0.25 ]{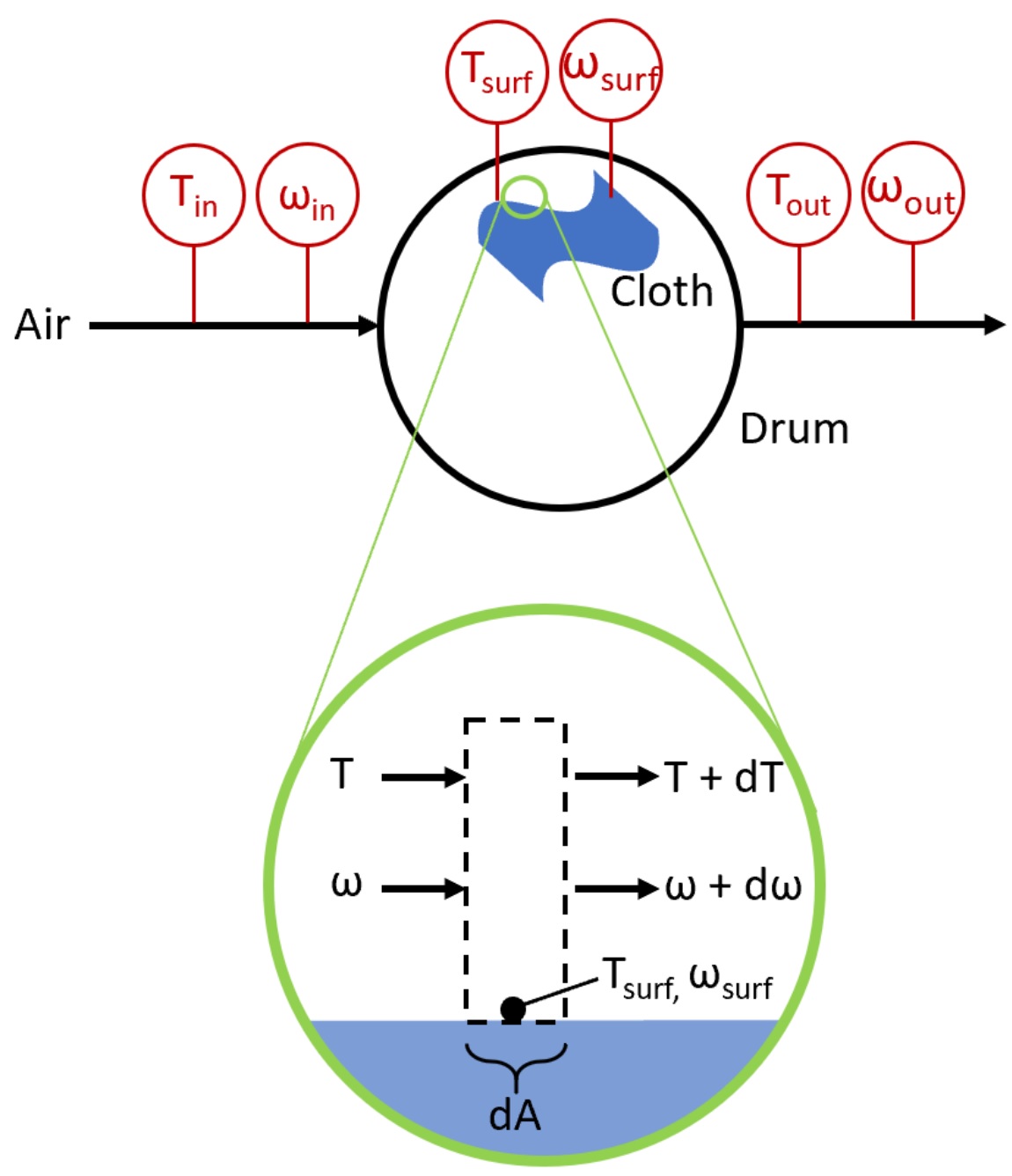}
\caption{Macro-level (top) and differential view (bottom) of the heat and mass transfer process in a dryer drum \citep{gluesenkamp2019efficient}}
\label{glupic2}
\end{figure}

The provided numerical methodology by Deans helped further researchers to model clothes drying cycles followingly \citep{stawreberg2013potential, somdalen2018theoretical, huang2019simulation}. Gluesenkamp utilized almost the same methodology to evaluate dying performance \citep{gluesenkamp2018experimental, gluesenkamp2019efficient}. As shown in Fig. \ref{glupic2}, for a piece of clothing in the drum the evolution of the sensible and latent energies can be described as follows:

\begin{equation}
\dot{m}_a c_p d T=h\; dA \left(T_{\text {surf }}-T\right) \\
\label{glueNew1}
\end{equation}

\begin{equation}
\dot{m}_a h_{f g} d \omega=k \; d A\left(\omega_{\text {surf }}-\omega\right) h_{f g} \\
\label{glueNew2}
\end{equation}

where $T_{\text {surf}}$ and $\omega_{\text {surf }}$ are the temperature and humidity ratio, respectively, of the cloth surface and are assumed to be uniform all over the drum. $A$ and $h_{f g}$ are the wet surface area and water latent heat of evaporation, respectively. For each of the ordinary differential equations, integration across the whole surface can be applied under the assumption that the surface temperature is spatially uniform \citep{gluesenkamp2019efficient}. This yields Eqs. \ref{glueNew3} and \ref{glueNew4} for heat and mass transfer, respectively:

\begin{equation}
T_{\text {out }}=T_{\text {surf }}-\left(T_{\text {surf }}-T_{\text {in }}\right) e^{-N T U} \\
\label{glueNew3}
\end{equation}

\begin{equation}
\omega_{\text {out }}=\omega_{\text {surf }}-\left(\omega_{\text {surf }}-\omega_{\text {in }}\right) e^{-N T U} \\
\label{glueNew4}
\end{equation}

where NTU is defined as the number of transfer units \citep{shen2016heat}:

\begin{equation}
N T U=\frac{h A_t}{\dot{m}_a c_p} \;,\\
\label{glueNew5}
\end{equation}

or as a general rule for the heat exchanger \cite{bergman2011introduction}:

\begin{equation}
\mathrm{NTU}=-\ln (1-\varepsilon) \\
\label{glueNew6}
\end{equation}

By rewriting Eqs. \ref{glueNew3} and \ref{glueNew4}, they could introduce two parameters namely heat and moisture transfer effectiveness, to address the efficiency of the drying cycle:

\begin{equation}
\varepsilon_{\mathrm{H}}=1-\frac{\mathrm{T}_{\text {surf }}-\mathrm{T}_{\text {out }}}{\mathrm{T}_{\text {surf }}-\mathrm{T}_{\mathrm{in}}} \\
\label{glueeqn1}
\end{equation}

\begin{equation}
\varepsilon_{\mathrm{M}}=1-\frac{\omega_{\text {surf }}-\omega_{\text {out }}}{\omega_{\text {surf }}-\omega_{\text {in }}} \\
\label{glueeqn2}
\end{equation}

Indices $in$, $out$ and $surf$ represent regions of inlet and outlet of drum and clothes surface. Temperature and specific humidity were measured directly with sensors at the inlet and outlet of the drum. As clothes temperature and humidity are hard to measure, they employed two following assumptions. 

First, the clothing surface is assumed to be saturated:

\begin{equation}
\omega_{\text {surf }}=\omega_{\text {surf,sat }} \\
\label{glueeqn3}
\end{equation} 

This allows the humidity ratio to be calculated as a function of surface temperature. Combining Eq. \ref{glueeqn2} and \ref{glueeqn3}:

\begin{equation}
\varepsilon_{\mathrm{M}}=1-\frac{\omega_{\text {surf,sat }}-\omega_{\text {out }}}{\omega_{\text {surf,sat }}-\omega_{\text {in }}} \\
\label{glueeqn4}
\end{equation}

$\omega_{\text {surf,sat }}$ is function of saturated vapor pressure at cloth's surface temperature (i.e. $P_{w,sat}=f(T_{surf})$) and atmospheric pressure, $P_{a}$, \citep{gluesenkamp2018experimental}:

\begin{equation}
\omega_{\text {surf,sat }}=\frac{\left.0.622 P_{w,sat}\right|_{T_{\text {surf }}}}{\left(P_{a}-\left.P_{w,sat}\right|_{T_{\text {surf }}}\right)} \\
\label{glueeqn5}
\end{equation}

Recall that Lewis number characterizes fluid flows when there is both convective heat and mass transfer:

\begin{equation}
\mathrm{Le}=\frac{\mathrm{Sc}}{\mathrm{Pr}}=\frac{\alpha}{D}=\frac{\text { Thermal diffusivity }}{\text { Mass diffusivity }} \\
\label{LeEq}
\end{equation}

where $D$ is the diffusion coefficient. Sc and Pr are Schmidt and Prandtl numbers, respectively. By assuming the value for the heat and mass transfer effectivenesses to be equal (assumption of Lewis number of unity), Eqs. \ref{glueeqn1} and \ref{glueeqn2} become equal as well:

\begin{equation}
\varepsilon_{\mathrm{H}}=\varepsilon_{\mathrm{M}}
\label{glueeqn6}
\end{equation}

The reduction of variables helped them estimate the efficiency of the drying cycle by solving Eqs. \ref{glueeqn1}-\ref{glueeqn6}. 

\begin{figure*}
\centering  \vspace*{0cm} \hspace*{0cm}
\includegraphics[scale=0.5 ]{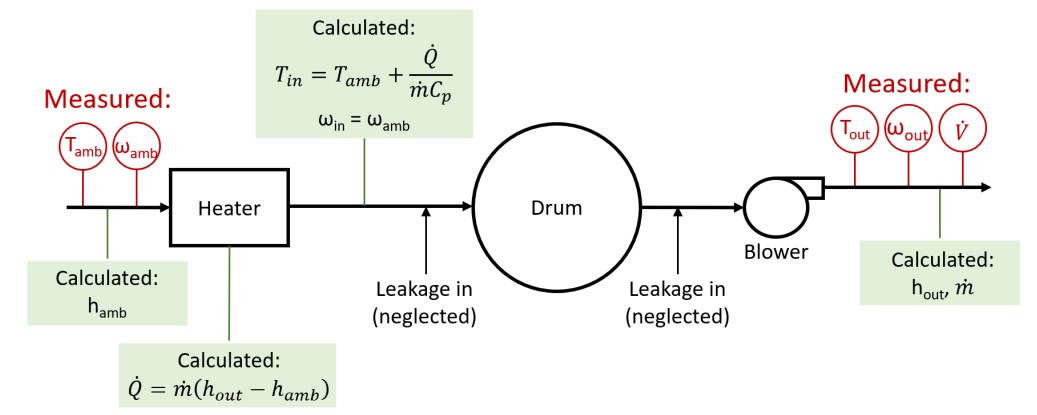}
\caption{Process and instrumentation diagram of experimental evaluation in \citep{gluesenkamp2018experimental}}
\label{glupic1}
\end{figure*}

\begin{figure}
\centering  \vspace*{0cm} \hspace*{0cm}
\includegraphics[scale=0.3 ]{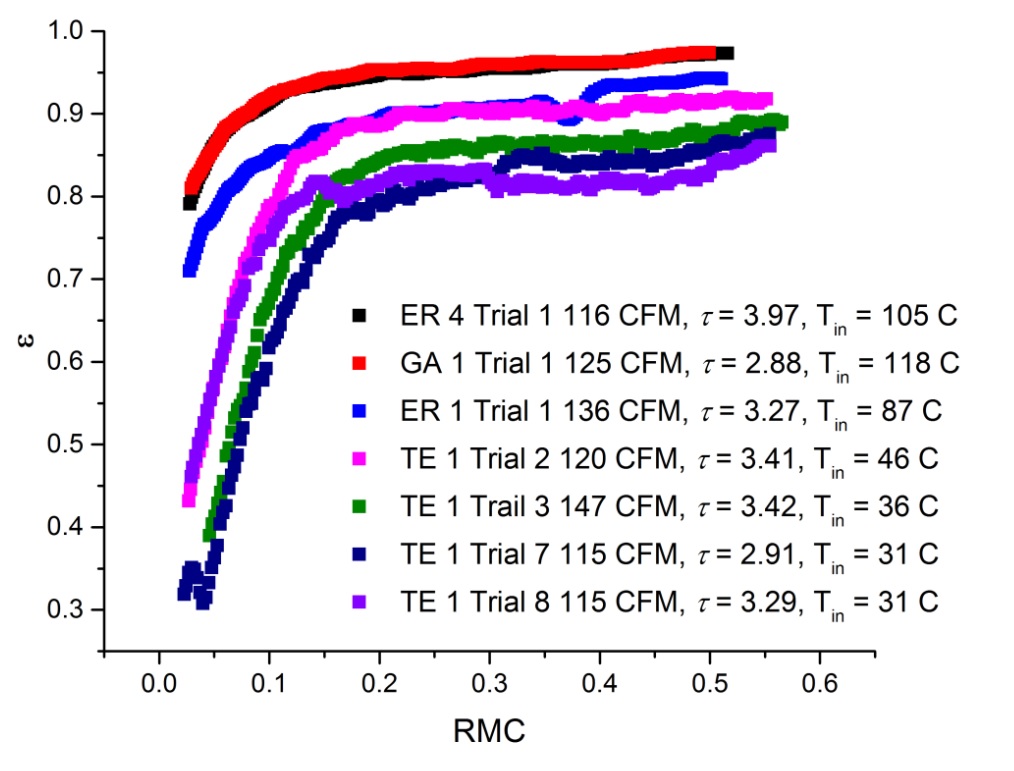}
\caption{Effectiveness in terms of remaining moisture content (RMC) at different drum residence times, fan flow rate and drum inlet temperature \citep{gluesenkamp2019efficient}}
\label{glupic3}
\end{figure}

Gluesenkamp estimated the effectiveness of the drying for different trials for clothes' dry load of 3kg, as shown in Fig. \ref{glupic3}. Accordingly, an increase in residence time has little effect on the effectiveness. Whilst, the airflow rate and inlet temperature increase the effectiveness through all stages of the drying.

Afterwards, they introduced a dimensional analysis for regression of the evaporation rate which is discussed in section \ref{Regression}.

The Deans's model was based on the 0-dimensional numerical analysis for solving energy balance. Some authors considered the temperature variation in the axial direction in the drum by solving conductive heat transfer between clothes \citep{wei2017mathematical}. However, the proposed 1-dimensional system of equations has shown no significant improvement in accuracy compared to the conventional 0-dimensional model.   

\begin{figure} 
\centering  \vspace*{0cm} \hspace*{0cm}
\includegraphics[scale=0.3]{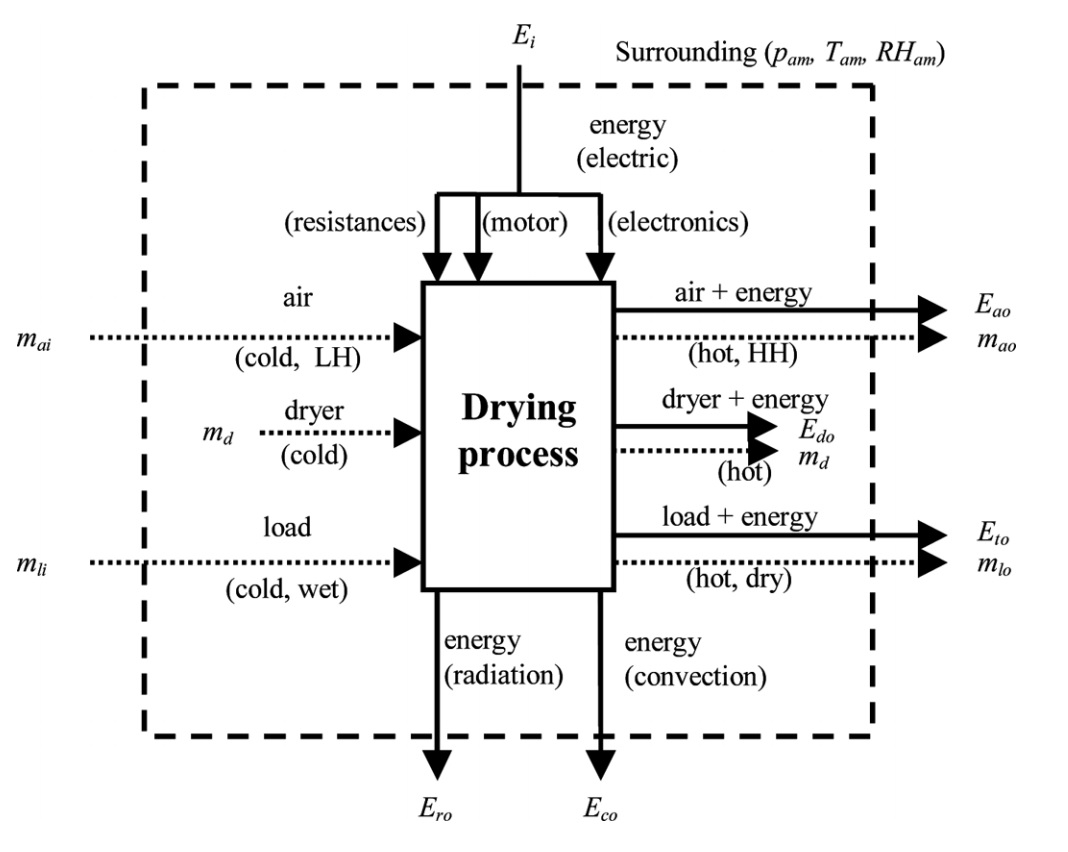}
\caption{Drying process scheme showing mass (dotted arrows) and energy fluxes (continuous arrows) \citep{huelsz2013total}}
\label{cvdryer}
\end{figure}

It is worth mentioning that plenty of studies conducted a similar approach to that of Deans and Lambert (i.e. energy and mass balance) to evaluate dryer performance. The difference lies in the details they include in the process. For instance, Conde evaluated dryers with a heat-recovery heat exchanger (typically cross-flow air-to-air heat exchanger for preheating the air stream at the inlet of the drum) \citep{conde1997energy}. Huelsz et al. included both radiative and convective heat losses of the dryer components to attain more accuracy \citep{huelsz2013total}. The schematic of their energy/mass balance model is shown in Fig. \ref{cvdryer}. Compared to that of Deans (\ref{deansnodes}), the difference is the neglect of the air bypass and the addition of different losses for the control volume. Jian and Zhao utilized a plate-fin heat exchanger for the condenser part of the dryer \citep{jian2017drying} and reported the effect of the water flow rate on the drying performance. Lee et al. investigated the drying performance in the presence of the heat pump \citep{lee2019rationally}. They assessed drying efficiency based on the following relation \citep{prasertsan1996heat} which is used frequently in related studies:

\begin{equation}
\eta_{d r y}=\frac{T_{i n}-T_{o u t}}{T_{\text {in }}-T_{o u t}}=\frac{\omega_{o u t}-\omega_{i n}}{\omega_{\text {sat }}-\omega_{\text {in }}} \\
\end{equation}

Another frequently reported parameter is the specific moisture extraction rate (SMER) \citep{lee2009investigation}:

\begin{equation}
\text { SMER }=\frac{\dot{m}_{ev}}{PO_{tot}}
\label{SMEREqu}
\end{equation}

\subsection{Yi Model}
\begin{figure*} 
\centering  \vspace*{0cm} \hspace*{0cm}
\includegraphics[scale=0.4]{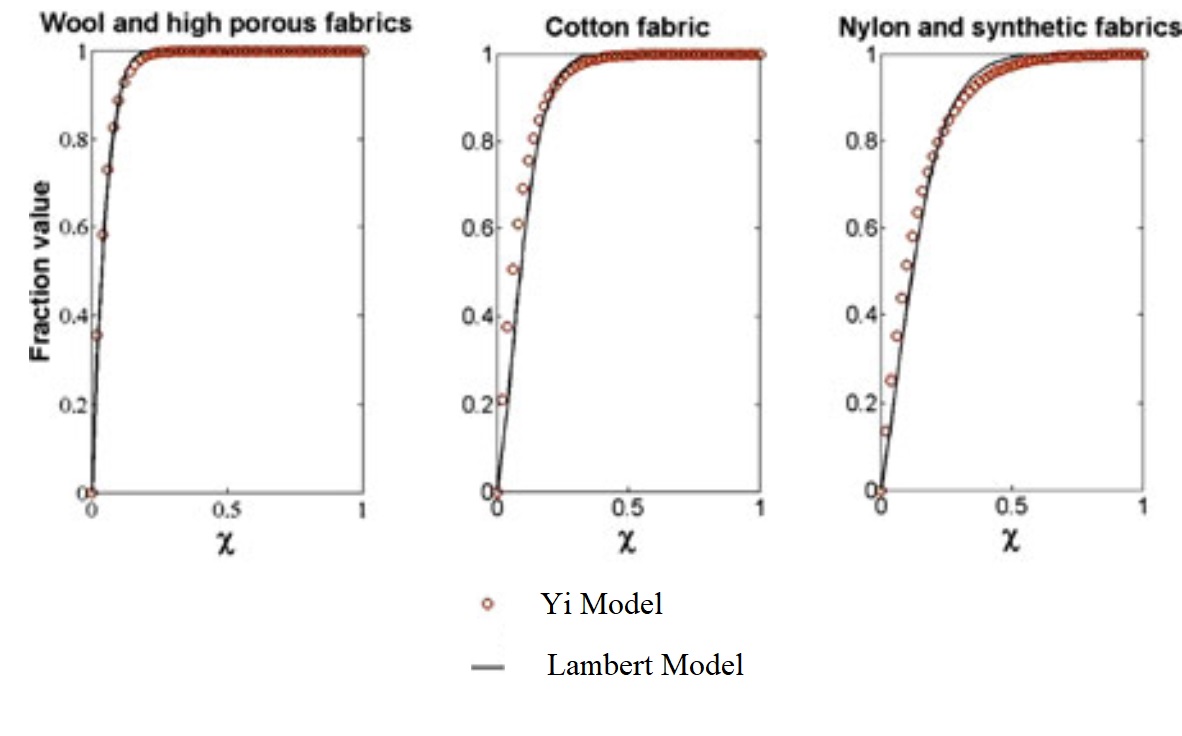}
\caption{The activity coefficient based on the numerical (Yi) and experimental (Lambert) methods for different fabrics \citep{yi2015new}}
\label{YiActivity}
\end{figure*}

Determining a correct expression for the specific humidity of clothes, $ \omega_{\mathrm{M}}$, in Eq. \ref{deansEq4} is of great importance. It is acceptable to consider it equal to the saturation-specific humidity, i.e. $ \omega_{\mathrm{M}} = \omega_{sat} $. Although previous approaches employed the conventional isotherm diagram (Fig. \ref{soriso2}) to account for the effect of the fabric's material, this method is based on the regression of experimental data (\citep{krischer2013wissenschaftlichen}). Yi, instead, suggested the utilization of probability density functions (PDF) based on the Boltzmann constant to explain underlying physics (\cite{yi2015new}). To this end, they initially quantified the amount of non-evaporated water in the clothes $m_{w, \mathrm{cl}}^{\text {no -evap }}$ as:

\begin{equation}
m_{w, \mathrm{cl}}^{\text {no -evap }}\left(\beta, \chi, T_{\mathrm{clo}}\right)=m_{w, \mathrm{clo}}^{\max } a\left(\beta, \chi, T_{\mathrm{clo}}\right)\\
\end{equation}

\begin{equation}
a\left(\beta, \chi, T_{\mathrm{clo}}\right)=1-\exp \left(-\frac{\beta \chi Q_{\mathrm{LH}}}{k_B T_{\mathrm{clo}}}\right)\\
\end{equation}

where $m_{w, \mathrm{clo}}^{\max }$ is the maximum amount of water in the clothes, $a\left(\beta, \chi, T_{\mathrm{clo}}\right)$ is the activity coefficient, $\beta$  is the property parameter to characterize the type of clothes fabrics, $Q_{\mathrm{LH}}$ is the latent heat of vaporization, and $k_B$ is the Boltzmann constant \citep{huang2008statistical}. Hence, they could estimate the water activity coefficient numerically. A comparison between the water activity coefficient for different clothes types is shown in Fig. \ref{YiActivity}. Clearly, their model was able to track the experimental trend.

\begin{figure} 
\centering  \vspace*{0cm} \hspace*{0cm}
\includegraphics[scale=0.35]{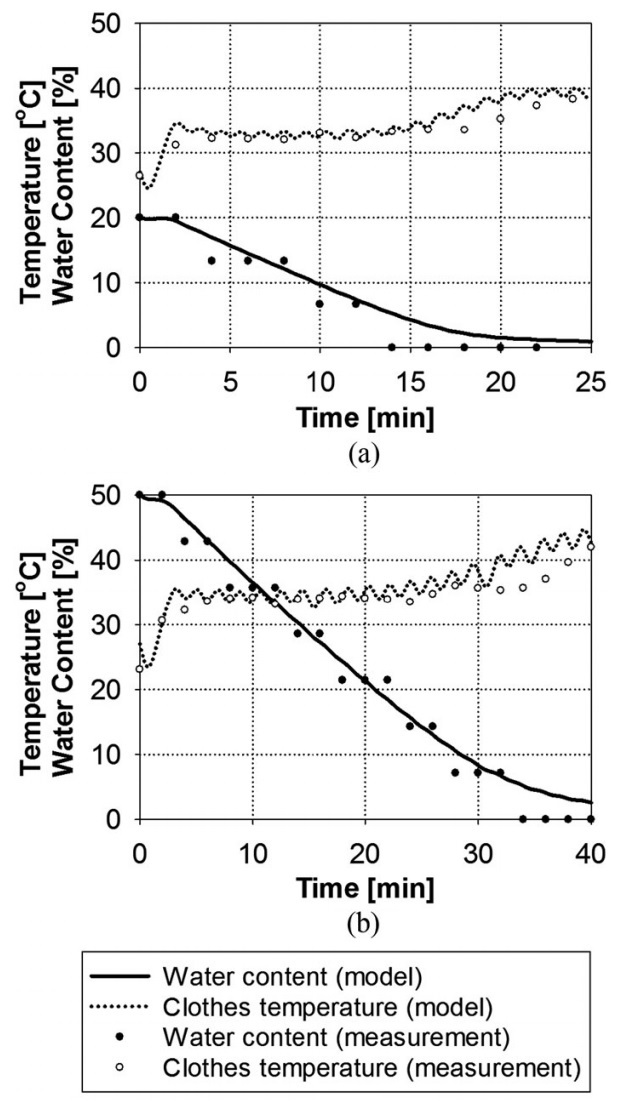}
\caption{Comparison of the clothes temperature and water content between the simulation and the experiments for (a) synthetic fabrics and (b) cotton fabrics \cite{yi2015new}}
\label{YiFig1}
\end{figure}

\begin{figure} 
\centering  \vspace*{0cm} \hspace*{0cm}
\includegraphics[scale=0.3]{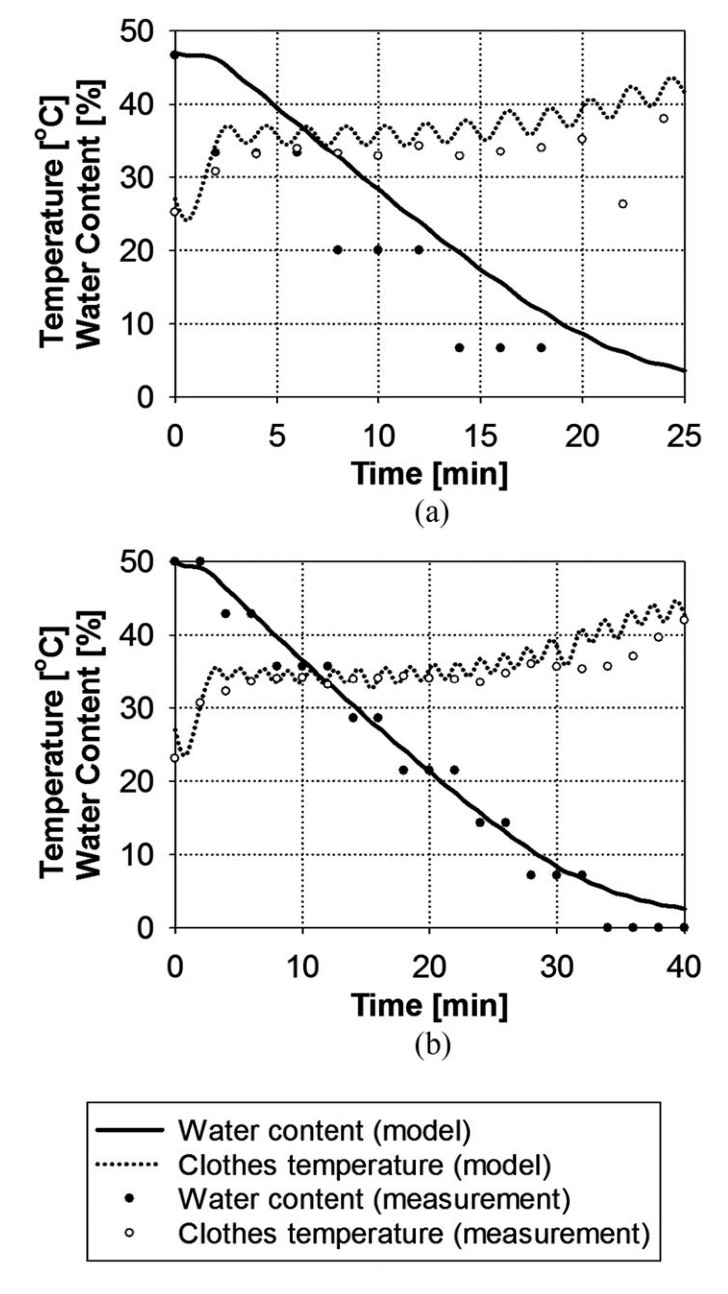}
\caption{Comparison of the clothes temperature and the water content of different clothes loads: (a) 1.5 kg and (b) 2.8 kg \cite{yi2015new}}
\label{YiFig2}
\end{figure}

Estimation of the evaporation rate in the drum was based on the fluid film analogy on the clothing surface. Accordingly, for the adjacent film, the amount of water in the air cannot exceed its saturated value. Thus, the $\omega_{\mathrm{M}}$ can be given by the following inequality:

\begin{equation}
\gamma m_{w, \mathrm{clo}}^{\max }\left(T_{\mathrm{clo}}\right) a\left(\beta, \chi, T_{\mathrm{clo}}\right)\\
\geq m_{\text {air }} \omega_{w, \text { air }}^{\text {sat }}\left(T_{\mathrm{clo}}\right) a(\beta, \chi, T)
\end{equation}

where $\gamma$ is the desorption weighting factor at the clothes interfaces and its value is less than 1. They solved the set of differential equations numerically using MATLAB software. They have validated their results for two types of clothing (Fig. \ref{YiFig1}) and different mass loads (Fig. \ref{YiFig2}). As shown, both synthetic and cotton fabrics show a similar range of temperatures through the drying. Moreover, the mass load has a linear effect on the drying time. 

\subsection{Regression}  \label{Regression}
All the previous methodologies made use of Chilton-Colburn analogy to investigate mass transfer in the drum. However, some authors tried to utilize the dimensional analysis based on Sherwood number to relate machine parameters to the evaporation rate. Bassily's work is considered to be one of the pioneer studies in this regard \citep{bassily2003correlation}. 

The Sherwood number, $\mathrm{Sh}$, (also called the mass transfer Nusselt number) is the ratio of the convective mass transfer to the rate of diffusive mass transport \citep{cochran2009condensing}. 

\begin{equation}
\mathrm{Sh}=\frac{h_{evap} L}{D}=\frac{\text { Convection mass transfer }}{\text { Diffusion rate }}\\
\end{equation}

where $h_{evap}$ and $L$ are convective mass transfer coefficient and length scale, respectively. 

Bassily developed a dimensionless correlation for estimating area-Sherwood number $\mathrm{Sh_{a}}$ (which will be defined later on) as a function of drum angular velocity, $\mathrm{N_{d}}$, textile weight, $\mathrm{M_{t}}$, Reynolds number, $\mathrm{Re}$, Schmidt number, $\mathrm{Sc}$, and Gukhman number, $\mathrm{Gu}$. The definition of Schmidt and Gukhman numbers are as follows \citep{smolsky1962heat, kovaci2022evaluation}:

\begin{equation}
\mathrm{Sc}=\frac{\nu}{D} = \frac{\text { Momentum diffusivity }}{\text { Diffusion rate }}\\
\end{equation}

\begin{equation}
\mathrm{Gu}=\frac{T_{g}-T_{gw}}{T_{g}} \\
\end{equation}

where $T_{g}$ and $T_{gw}$ are gas dry and wet bulb temperatures, respectively.

They found the following relation for the $\mathrm{Sh}_{a}$ using 32 experimental runs:

\begin{equation}
\mathrm{Sh}_{a} = 0.028 M_{t}^{0.58} N_{d}^{0.06} \mathrm{Gu}^{-1.38} \mathrm{Re}^{0.35} \mathrm{Sc}^{1.22}
\label{bassilyeq2}
\end{equation}

Then area–Sherwood number $\mathrm{Sh}_{a}$ can be used to determine area-mass transfer coefficient, $kA$, as follows \citep{rasti2021review}:

\begin{equation}
k A =\frac{\mathrm{Sh}_{a} \times D}{d_{d}} \;,
\label{bassilyeq1}
\end{equation}

where $d_{d}$ and $D$ are drum diameter and diffusion coefficient, respectively. 

The diffusion coefficient is a function of the instantaneous moisture content, $C(t)$:

\begin{equation}
D =-1.71778 \times 10^{-7}\left(\frac{C(t)}{1.2727}\right)^{5}
\label{deanEqDab}
\end{equation}

\begin{figure*} 
\centering  \vspace*{0cm} \hspace*{0cm}
\includegraphics[scale=0.4]{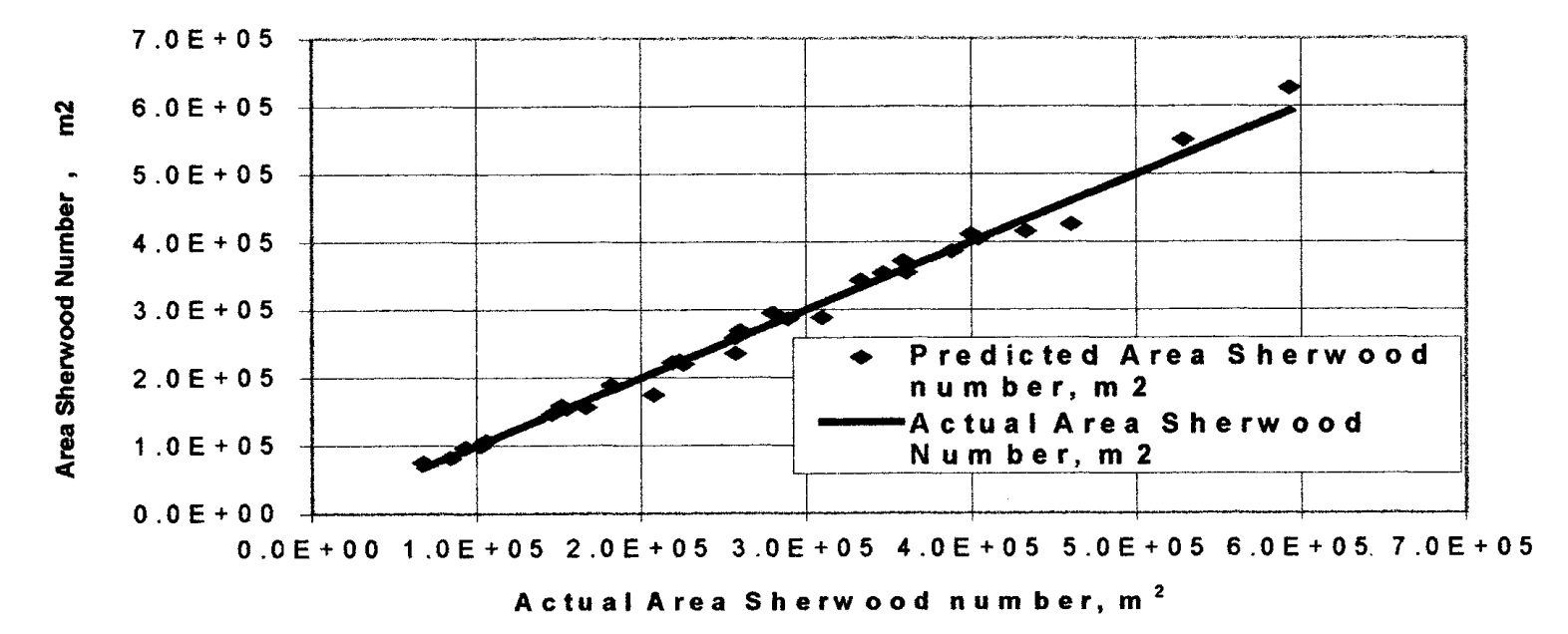}
\caption{The measure and  predicted area Sherwood numbers \citep{bassily2003correlation}}
\label{Bassily1}
\end{figure*}

The comparison between the measured and correlation area-Sherwood numbers are shown in Fig. \ref{Bassily1}.

More recent literature tends more and more to regress the drying indices by the machine inputs \citep{gluesenkamp2019efficient, novak2019influence}. Noval et al. introduced relationships for drying time, evaporation rate and area-mass transfer coefficients in terms of volumetric flow rate, inlet temperature and relative humidity as follows:

\begin{equation}
t=39434 \times \dot{V}_{\mathrm{in}}^{-0.9065} T_{\mathrm{in}}^{-0.8338} \mathrm{RH}_{\mathrm{in}}^{0.1598} \\
\label{novaleq1}
\end{equation}

\begin{equation}
\dot{m}_{\text {evap }} =5.53 \times 10^{-7} \dot{V}_{\mathrm{in}}^{0.8305} T_{\mathrm{in}}^{0.8066} \mathrm{RH}_{\mathrm{in}}^{-0.1914}\\
\label{novaleq2}
\end{equation}

\begin{equation}
k A =9.06 \times 10^{-8} \dot{V}_{\mathrm{in}}^{0.6824} T_{\mathrm{in}}^{-0.5561} \mathrm{RH}_{\mathrm{in}}^{-0.0293}
\label{novaleq3}
\end{equation}

Gluesenkamp derived a correlation for predicting the effectiveness of heat and mass transfer phenomena in the drum. The effectiveness is a function of five parameters: the textile weight $m_{c}$, the water content in textiles $m_{w}$, the air weight in the drum $m_{a}$, the air inlet flow rate $\dot{m}_a$, and the falling time (which is the time taken for a piece of clothing to fall from the top to the bottom of the drum) $t$:

\begin{equation}
\varepsilon=f\left(m_{clo}, m_w, m_{a}, \dot{m}_a, t\right) \\
\label{gluenew2-1}
\end{equation}

By means of Buckingham $\pi$-theorem they could find the following dimensionless parameters:

\begin{equation}
\pi_1=\varepsilon \\
\pi_2=\frac{\dot{m}_a t}{m_{a}} \\
\pi_3=\frac{m_{clo}}{m_{a}} \\
\pi_4=\frac{m_{w}}{m_{a}} \\
\label{gluenew2-2}
\end{equation}

They took the first-order polynomial for $\pi_2$ and $\pi_3$ where a forth-order polynomial is considered for $\pi_4$:

\begin{equation}
\pi_1=a+b_1 \pi_2+c_1 \pi_3+d_1 \pi_4+d_2 \pi_4^2+d_3 \pi_4^3+d_4 \pi_4^4
\label{gluenew2-3}
\end{equation}

They did a set of experiments to find the coefficients of Eq. \ref{gluenew2-3} for HLD-1-1992 (100\% cotton-like sheets, table cloths, bath towels, long sleeve shirts, T-shirts, pillowcases, shorts, washcloths, handkerchiefs \citep{AHAM1}) and  HLD-1-2009 (100\% cotton including sheets, pillowcases, towels \citep{AHAM2}). The values of coefficients are shown in Table \ref{tab:gluenew2-1}.

\begin{table}[]
\centering
\renewcommand{\arraystretch}{2}
\caption{Coefficients in Eq. \ref{gluenew2-3} \citep{gluesenkamp2019efficient}}
\label{tab:gluenew2-1}
\begin{tabular}{|c|c|c|}
\hline
\textbf{Coefficient} & \textbf{AHAM 2009} & \textbf{AHAM 1992} \\ \hline
\textbf{$a$}           & 0.677584           & 0.354049           \\ \hline
\textbf{$b_1$}          & -4.456546          & -2.11532           \\ \hline
\textbf{$c_1$}          & 0.015385           & 0.00212            \\ \hline
\textbf{$d_1$}          & 0.290916           & 0.403633           \\ \hline
\textbf{$d_2$}          & -0.071735          & -0.087495          \\ \hline
\textbf{$d_3$}          & 0.007304           & 0.008034           \\ \hline
\textbf{$d_4$}          & -0.000249          & -0.000252          \\ \hline
\end{tabular}
\end{table}

There has been no report on the sensitivity analysis for the dryer's input parameter, and therefore, there is no agreement on the independent variables in the regression models. For instance, Gluesenkamp did not take drum speed or air temperature into account \citep{gluesenkamp2019efficient}.

 In addition, although the implementation of the regressed equations is fast as it does not require solving any differential equation, the validity of the correlations holds only for the reference machine. Put differently, the mentioned equations cannot be generalized for other tumble dryers. In this matter, Rasti and Jeong compared the mentioned regression models of Bassily, Gluesenkamp and Novak (\citep{rasti2021review}). For the lowest and medium clothes load the Gluesenkamp and Novak models, respectively, could capture the physics of evaporation. However, at the high load (9kg), none of the models could anticipate the phenomena. Worth mentioning that the Bassily model failed to capture the physics in all the cases for that specific machine.

\begin{figure} 
\centering  \vspace*{0cm} \hspace*{0cm}
\includegraphics[scale=0.38]{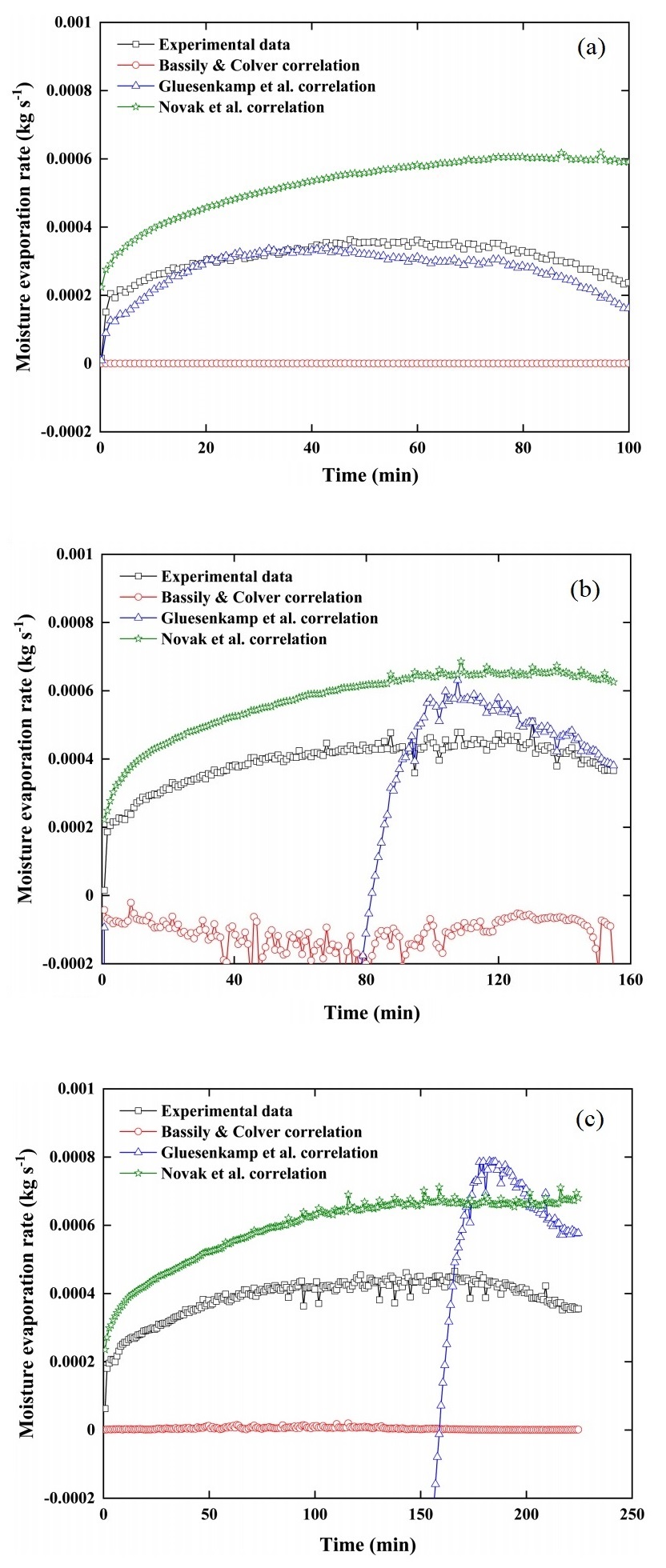}
\caption{A comparison between the different regression models for the evaporation rate for different loads; (a) 3kg (b) 6kg and (c) 9kg \citep{rasti2021review}}
\label{RegressionModels}
\end{figure}

\subsection{Dynamics of Motions}

\begin{figure*} 
\centering  \vspace*{0cm} \hspace*{0cm}
\includegraphics[scale=0.4]{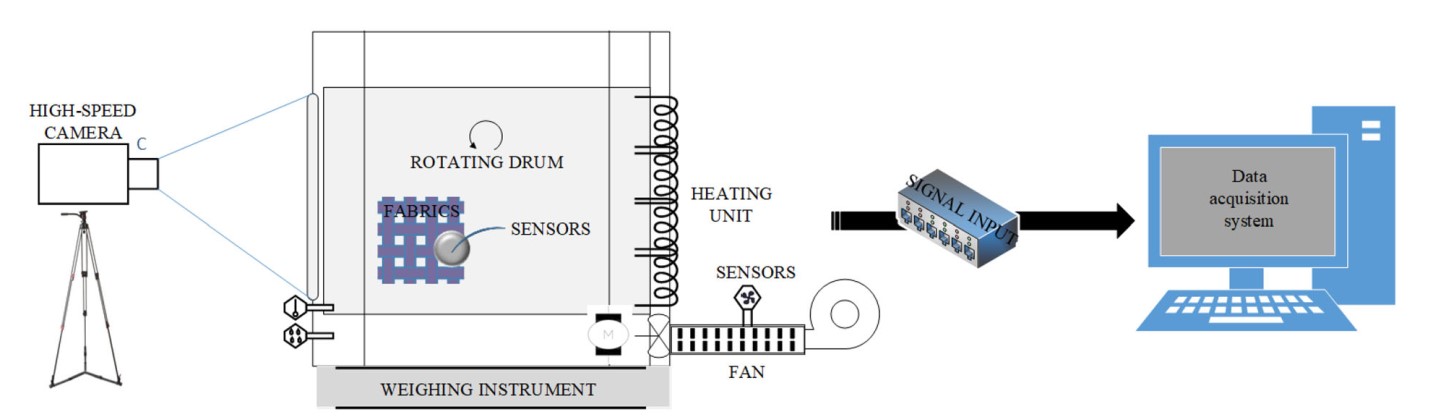}
 \caption{Schematic of the video capturing and processing system \citep{yu2021dynamics}}
\label{ImagingTech}
\end{figure*}

With the advances in measuring technologies, a newer branch of study in the field of tumble dryers has arisen in recent years. In this regard, some authors tried to relate the thermodynamics of evaporation to the dynamics of the clothes in the drum. Modern optical methods and high-speed video capturing instruments, such as Positron Emission Particle Tracking (PEPT), for flow and clothes visualization can be utilized to track fabrics during the drying process. Then, with the aid of image processing techniques, mainly Neural Network (NN) methods the gravity centre of the tracer fabrics and their motion characteristics can be derived. The typical procedure is shown in Fig. \ref{ImagingTech}. 

\begin{figure*} 
\centering  \vspace*{0cm} \hspace*{0cm}
\includegraphics[scale=0.45]{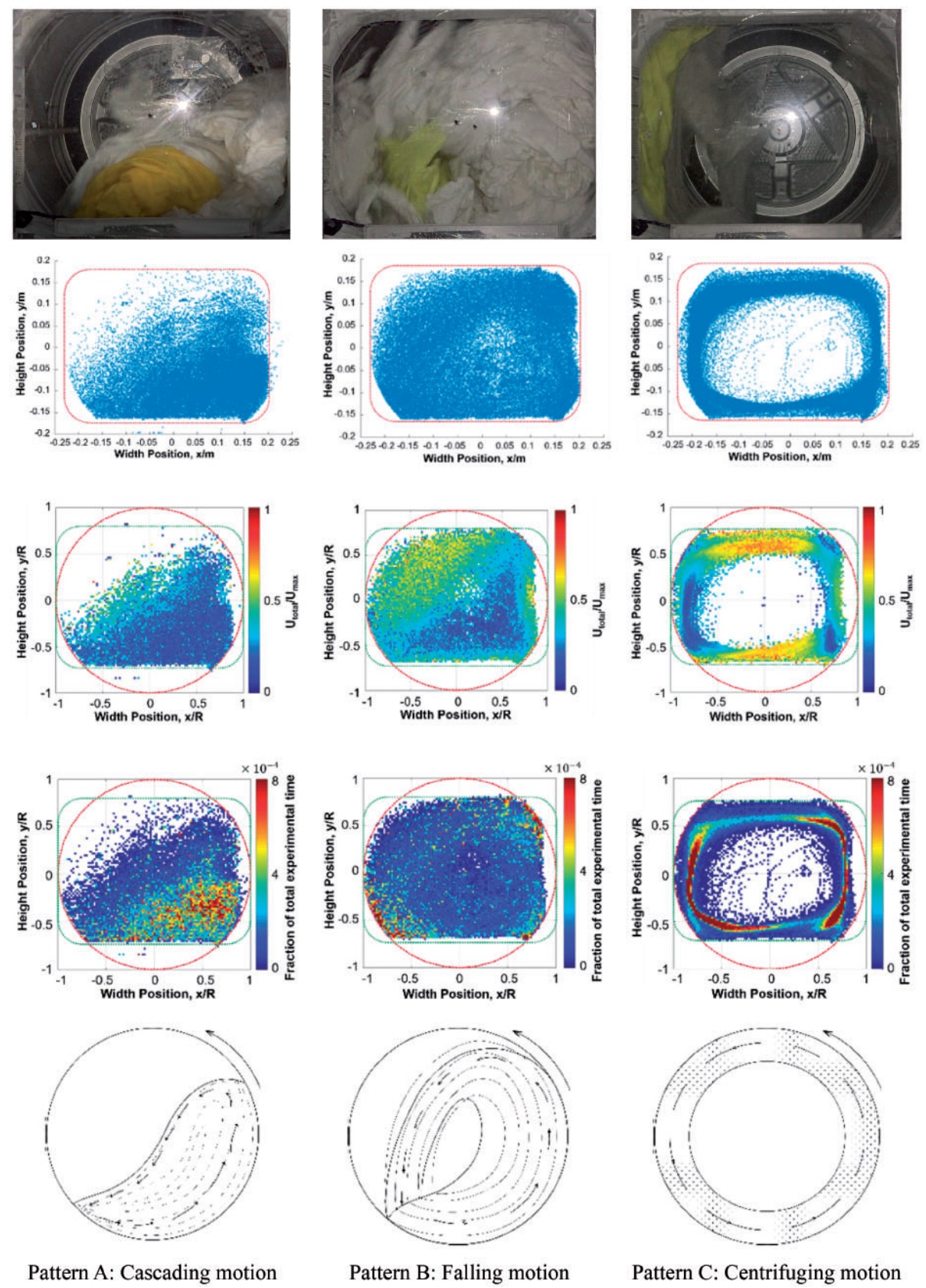}
 \caption{Three categories of cotton textile transverse motion. From top to bottom: experimental image, trajectories, velocity distribution, residence time distribution and schematic of each motion pattern \citep{yu2021dynamics}}
\label{cotimag}
\end{figure*}

Yu et al. derived relations for the drying rate and SMER in terms of clothes motion \citep{yu2021dynamics}. They reported SMER as follows:

\begin{equation}
\text{SMER} = -4.1 -0.1 x_{1} + 0.02x_{2} + 0.69x_{3} \;,
\label{yuimaproceq1}
\end{equation}

where $x_{1}$ and $x_{2}$ are the ratio of the total motion area of the tracer and the fabric spread area, respectively, and $x_{3}$ is the distance between the drum center and tracer fabric center. As shown in Fig. \ref{cotimag}, dynamics of clothes during the drying process including clothes trajectory, velocity and residence time are derived by image processing techniques. Consequently, additional results such as motion regime (based on what is explained in Fig. \ref{ImagingTechFroude}) are derived. Similar studies have been done by \citep{yu2021effects} and \citep{lee2022effect}.

This technique is also considered to be machine-specific. As there are only a few research studies relating to this technique, it is hard to assess it correctly compared to other methods.

\section{Evaporation Rate}
The mentioned methodologies include a major part of the literature regarding tumble dryers mainly focused on the prediction of the dryer working process. Interestingly, all the numerical and most of the experimental studies pursue a common goal, which is the evaluation of the \textit{evaporation rate}. With the knowledge of the evaporation rate, the important parameters for evaluation of the system performance can be determined. Among these parameters, SMER and drying time are the most important ones. Other auxiliary variables such as the drying evenness can be related to the evaporation rate \citep{yu2018wrinkling, wei2018enhancing, hu2016effects, chen2018comparisons}. Moreover, the drying rate has an effect on the clothes wrinkling during the drying cycle \citep{el2022state}.

Briefly, there exist four main methods for predicting this value in the literature. In the following, details of these methods will be presented.

\vspace{0.25cm}

\subsection{Direct Method}
\begin{figure}  
\centering  \vspace*{0cm} \hspace*{0cm}
\includegraphics[scale=0.25]{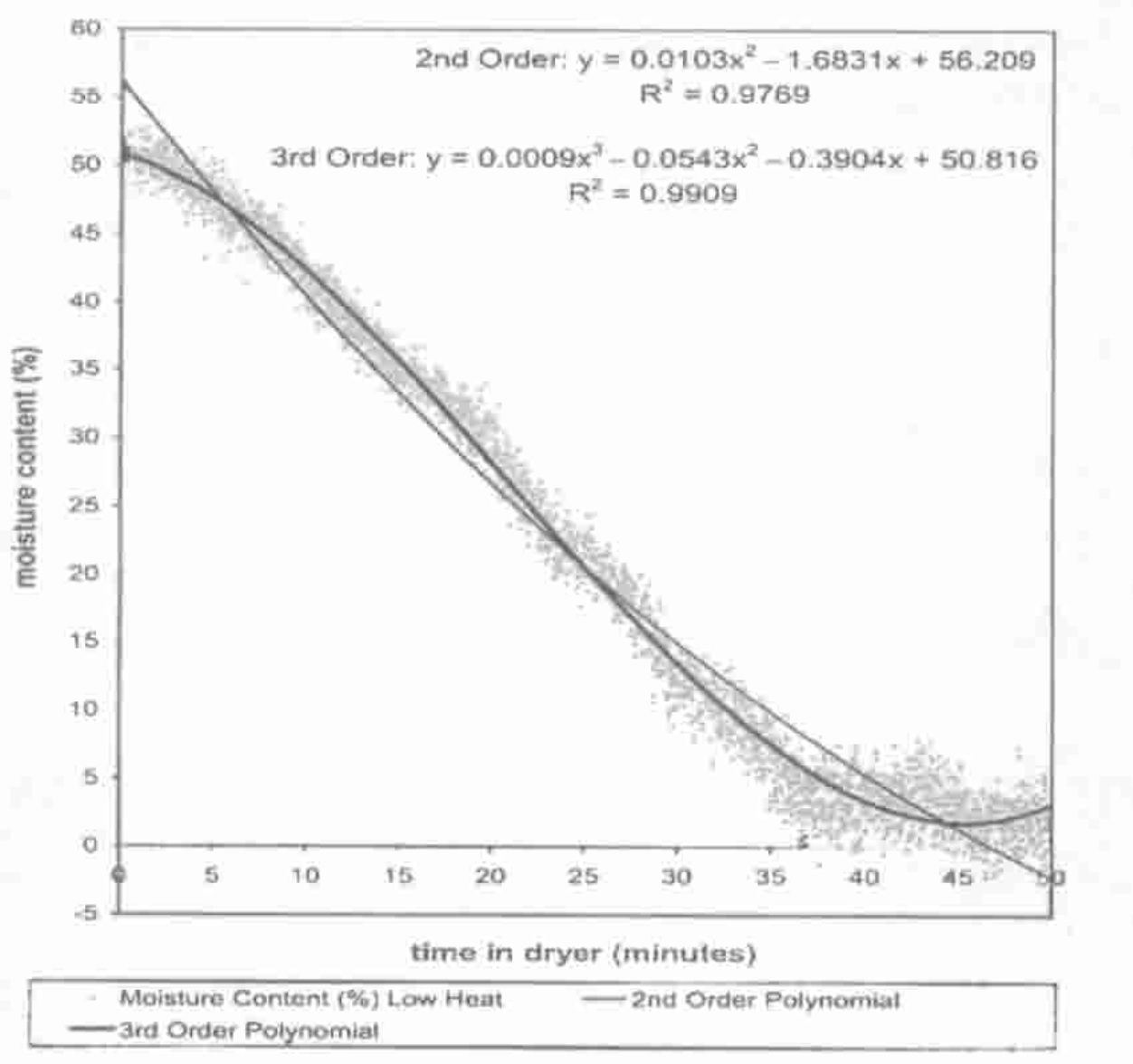}
\caption{Change in moisture content during tumble-drying at low-temperature setting \citep{higgins2003effect}}
\label{higgings1}
\end{figure}

This method aims to measure the weight of the components in the drum directly. As the drum rotates, the content experiences centrifugal and Coriolis forces continuously. Thus, high oscillations may exist in the measurement. This method is not popular in the literature, however, was utilized by some authors \citep{higgins2003effect}. Higgins et al. employed a platform scale with a sampling frequency of $1\si{sec}$ to continuously measure the drum contents. Fig. \ref{higgings1} shows the change of moisture content through time. As depicted, the resulting data is noisy but can still be used by inserting specific filters and averaging functions.

\vspace{0.15cm}

\subsection{Constant Rate}
There exist some reports in the literature that employed constant evaporation rate through the whole drying period \citep{bansal2001improving}. We et al. took an evaporation rate of $3.5 \frac{kg}{h}$ in their numerical method \citep{braun2002energy}.  

\begin{figure}  
\centering  \vspace*{0cm} \hspace*{0cm}
\includegraphics[scale=0.5]{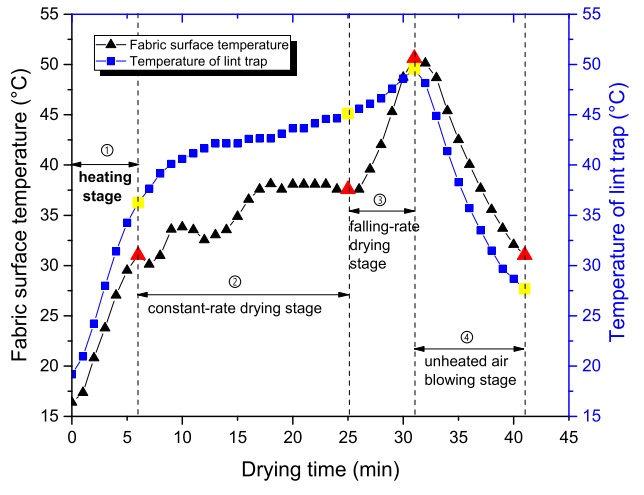}
\caption{The changes of drum temperature during the drying process \citep{yu2021effects}}
\label{yuFigTempChange}
\end{figure}

\begin{figure}  
\centering  \vspace*{0cm} \hspace*{0cm}
\includegraphics[scale=0.5]{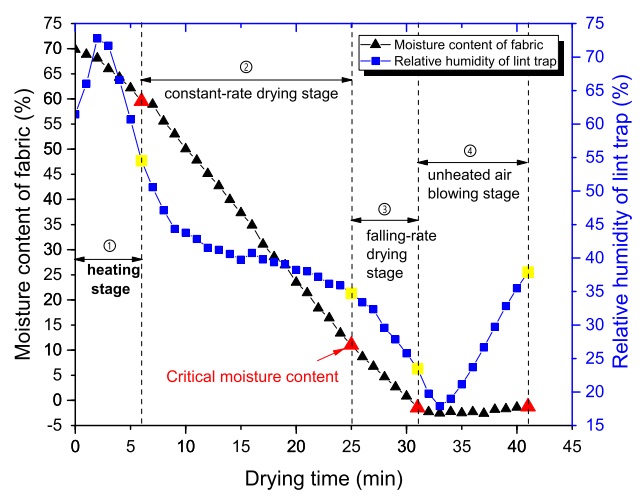}
\caption{The changes of moisture during the drying process \citep{yu2021effects}}
\label{yuFigMoisChange}
\end{figure}

The changes in temperature and moisture content in the drum have a non-linear trend through time. For instance, the typical evolution of the clothes temperature and moisture over time are depicted in Fig. \ref{yuFigTempChange} and \ref{yuFigMoisChange}. Implementation of the so-called constant rate for evaporation may cause errors in the estimation of the start (transient) and final (falling rate) stages of the drying. 

As mentioned before, the idea of a constant evaporation rate has been only used in studies accompanied by the heat pump cycle. In those situations, the main focus is on the heat pump's performance. Additionally, because of the complexity introduced to the numerical model by considering the heat pump, this assumption is acceptable. Thus, this method is not recommended for condenser tumble dryers.

\begin{figure}  
\centering  \vspace*{0cm} \hspace*{0cm}
\includegraphics[scale=0.45]{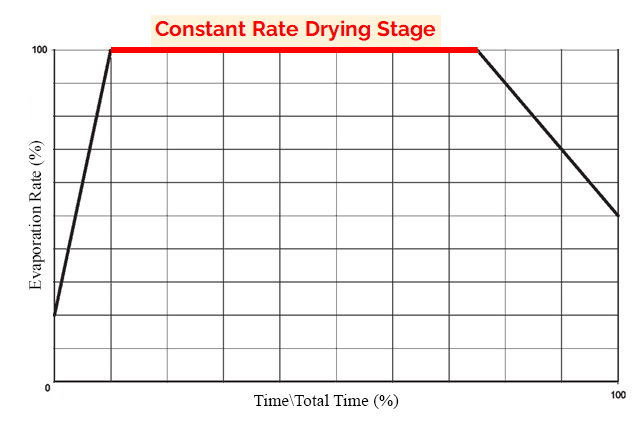}
\caption{Typical evaporation rate evolution in time}
\label{EvaRate3}
\end{figure}

\vspace{0.15cm}

\subsection{Chilton–Colburn Analogy}
Lambert as one of the pioneers in proposing a mathematical model for the tumble dryers estimated the evaporation rate (Eq. \ref{equlam}), where the area-mass transfer coefficient, $kA$, can be taken as a constant value, determined experimentally (regression) or by Sherwood analogy. The complete description has been provided previously.

\vspace{0.3cm}

\begin{figure*}  
\centering  \vspace*{0cm} \hspace*{0cm}
\includegraphics[scale=0.5]{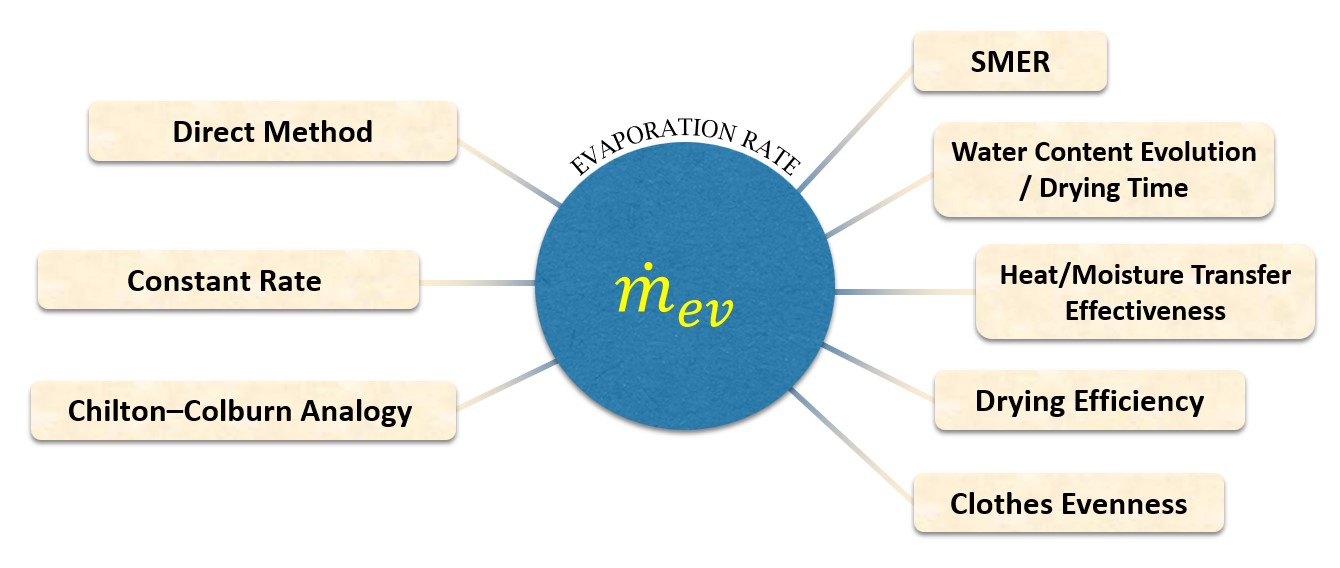}
\caption{Methods for evaluating the evaporation rate of the tumble dryer}
\label{EvapRateMeth}
\end{figure*}

The summary of mentioned methods for the evaporation rate and important output indices are shown in Fig. \ref{EvapRateMeth}. In Table. \ref{tableEvapModels} methodologies for evaluating evaporation rate for some other studies are presented as well.

\begin{table*}  
\caption{List of some employed methodologies for evaporation rate in previous studies}
\label{tableEvapModels}
\renewcommand{\arraystretch}{2}
\begin{tabular}{|c|c|c|}
\hline
\textbf{Researcher}    & \textbf{Year} & \textbf{Model}                        \\ \hline
Lee et al. \citep{lee2019rationally}             & 2019          & Lambert                               \\ \hline
Huang et al. \citep{huang2019simulation}           & 2019          & Lambert                               \\ \hline
Sian and Wang \citep{sian2019comparative}         & 2019          & Lambert                               \\ \hline
Yadav and Moon \citep{yadav2008modelling,yadav2008fabric}         & 2008          & Lambert                               \\ \hline
Deans \citep{deans2001modelling}                 & 2001          & Lambert                               \\ \hline
Yi et al. \citep{yi2015new}          & 2015          & Lambert with Modifications            \\ \hline
Cengel and Ghajar \citep{ccengel2020heat}               & 2014          & Lambert with Modifications            \\ \hline
Bengtsson et al. \citep{bengtsson2014performance}       & 2014          & Constant Rate                         \\ \hline
Stawreberg \citep{stawreberg2011energy}            & 2011          & Constant Rate                         \\ \hline
Stawreberg and Nilsson \citep{stawreberg2013potential} & 2013          & Constant Rate                         \\ \hline
Cranston et al. \citep{cranston2019efficient}       & 2019          & Constant Rate                         \\ \hline
Novak et al. \citep{novak2019influence}           & 2019          & Constant Rate                         \\ \hline
Higgins et al. \citep{higgins2003effect}           & 2003          & Direct Method                         \\ \hline
\end{tabular}
\end{table*}

\vspace{0.5cm}

\section{Conclusion}
Tumble dryers are now considered as common home appliances as they offer a convenient way of drying textiles. A critical review of the published literature concerning the analysis of domestic tumble dryers is presented. There has been a concentration during the decades on modelling heat and mass transfer through the dryer. In this regard, researchers have used various methodologies to predict the performance indices, particularly, the evaporation rate. 

The first mathematical model to determine the mass transfer in the drum was reported by Lambert in 1991. By utilizing Chilton–Colburn analogy they were able to define an area-mass transfer coefficient in terms of clothes' properties and the drum's inlet temperature and humidity to address the moisture evaporation rate. The evaporation rate can be considered to be the primary drying index with which other the drying time, efficiency, effectiveness and SMER can be derived, respectively. Moreover, Lambert divided evaporation into three stages namely transient, constant rate and falling rate through the drying cycle. 

During the last three decades, several experimental or numerical studies were published based on this concept, and furthermore, the model was then developed and parameterized by Deans, Gluesenkamp and Yi by considering more complexities and fewer assumptions. Furthermore, some literature focused on utilizing dimensional analysis and similarity to correlate the drying rate to the machine parameters. Bassily, for instance, developed a correlation for the area-Sherwood number as a function of clothes' mass, drum angular velocity, the mass flow rate of the dry air and geometrical parameters. More precise models were introduced by Novak and Gluesenkamp later on. The validity of the regressed has been presented as part of the review, and it was shown that this method lacks the required accuracy.

\pagebreak

\section{Bibliography}

\bibliographystyle{cas-model2-names}

\bibliography{ReferencesV1}


\end{document}